\documentclass{article}
\usepackage{amsmath}
\usepackage{amssymb}
\begin{document}
\def\e#1\e{\begin{equation}#1\end{equation}}
\def\ea#1\ea{\begin{align}#1\end{align}}
\def\eq#1{{\rm(\ref{#1})}}
\newtheorem{thm}{Theorem}[section]
\newtheorem{lem}[thm]{Lemma}
\newtheorem{prop}[thm]{Proposition}
\newtheorem{conj}[thm]{Conjecture}
\newtheorem{cor}[thm]{Corollary}
\newenvironment{dfn}{\medskip\refstepcounter{thm}
\noindent{\bf Definition \thesection.\arabic{thm}\ }}{\medskip}
\newenvironment{ex}{\medskip\refstepcounter{thm}
\noindent{\bf Example \thesection.\arabic{thm}\ }}{\medskip}
\newenvironment{proof}[1][,]{\medskip\ifcat,#1
\noindent{\it Proof.\ }\else\noindent{\it Proof of #1.\ }\fi}
{\relax\unskip\nobreak ~\hfill$\square$\medskip}
\def\dim{\mathop{\rm dim}}
\def\Re{\mathop{\rm Re}}
\def\Im{\mathop{\rm Im}}
\def\Ker{\mathop{\rm Ker}}
\def\Coker{\mathop{\rm Coker}}
\def\ind{\mathop{\rm ind}}
\def\sind{{\ts\mathop{\text{\rm s-ind}}}}
\def\lind{\mathop{\text{\rm l-ind}}}
\def\sign{\mathop{\rm sign}}
\def\End{\mathop{\rm End}}
\def\Vect{\mathop{\rm Vect}}
\def\vol{\mathop{\rm vol}}
\def\Hol{{\textstyle\mathop{\rm Hol}}}
\def\hol{{\textstyle\mathop{\mathfrak{hol}}}}
\def\Ric{\mathop{\rm Ric}}
\def\Vol{\mathop{\rm Vol}}
\def\Hess{\mathop{\rm Hess}}
\def\Image{\mathop{\rm Image}}
\def\Tr{\mathop{\rm Tr}}
\def\Spin{\mathop{\rm Spin}}
\def\Sp{\mathop{\rm Sp}}
\def\GL{\mathop{\rm GL}}
\def\SO{\mathop{\rm SO}}
\def\O{\mathop{\rm O}}
\def\U{\mathbin{\rm U}}
\def\SL{\mathop{\rm SL}}
\def\SU{\mathop{\rm SU}}
\def\ge{\geqslant} 
\def\le{\leqslant} 
\def\cal{\mathcal}
\def\H{\mathbb H\mskip1mu}
\def\R{\mathbb R\mskip1mu}
\def\Z{\mathbb Z\mskip1mu}
\def\N{\mathbb N\mskip1mu}
\def\C{\mathbb C\mskip1mu}
\def\sst{\scriptscriptstyle}
\def\sSi{{\smash{\sst\Si}}}
\def\sSii{{\smash{\sst\Si_i}}}
\def\sX{{\smash{\sst X}}}
\def\sXp{{\smash{\sst X'}}}
\def\sL{{\smash{\sst L}}}
\def\sN{{\smash{\sst N}}}
\def\D{{\cal D}}
\def\I{{\cal I}}
\def\M{{\cal M}}
\def\oM{{\,\,\ov{\!\!{\cal M}\!}\,}}
\def\O{{\cal O}}
\def\g{{\mathfrak g}}
\def\h{{\mathfrak h}}
\def\CP{\mathbb{CP}}
\def\u{\mathfrak{u}} 
\def\su{\mathfrak{su}} 
\def\so{\mathfrak{so}} 
\def\al{\alpha}
\def\be{\beta}
\def\ga{\gamma}
\def\de{\delta}
\def\ep{\epsilon}
\def\th{\theta}
\def\la{\lambda}
\def\ka{\kappa}
\def\vp{\varphi}
\def\si{\sigma}
\def\up{\upsilon}
\def\Up{\Upsilon}
\def\De{\Delta}
\def\La{\Lambda}
\def\Th{\Theta}
\def\Om{\Omega}
\def\Ga{\Gamma}
\def\Si{\Sigma}
\def\om{\omega}
\def\d{{\rm d}}
\def\pd{\partial}
\def\db{{\bar\partial}}
\def\ts{\textstyle}
\def\pha{\phantom}
\def\w{\wedge}
\def\lt{\ltimes}
\def\ti{\tilde}
\def\sm{\setminus}
\def\ov{\overline}
\def\ot{\otimes}
\def\bigot{\bigotimes}
\def\iy{\infty}
\def\ra{\rightarrow}
\def\t{\times}
\def\ha{{\textstyle\frac{1}{2}}}
\def\op{\oplus}
\def\bs{\boldsymbol}
\def\ms#1{\vert#1\vert^2}
\def\bms#1{\bigl\vert#1\bigr\vert^2}
\def\md#1{\vert #1 \vert}
\def\bmd#1{\bigl\vert #1 \bigr\vert}
\def\an#1{\langle#1\rangle}
\title{Lectures on special Lagrangian geometry}
\author{Dominic Joyce \\ Lincoln College, Oxford, OX1 3DR \\
{\tt dominic.joyce@lincoln.ox.ac.uk}}
\date{}
\maketitle

\section{Introduction}
\label{bp1}

{\it Calabi--Yau $m$-folds} $(M,J,\om,\Om)$ are compact complex
manifolds $(M,J)$ of complex dimension $m$, equipped with a
Ricci-flat K\"ahler metric $g$ with K\"ahler form $\om$, and a
holomorphic $(m,0)$-form $\Om$ of constant length $\ms{\Om}=2^m$.
Using Algebraic Geometry and Yau's solution of the Calabi
Conjecture, one can construct them in huge numbers. String
Theorists (a species of theoretical physicist) are very
interested in Calabi--Yau 3-folds, and have made some
extraordinary conjectures about them, in the subject
known as {\it Mirror Symmetry}.

{\it Special Lagrangian submanifolds}, or {\it SL $m$-folds}, are
a distinguished class of real $m$-dimensional minimal submanifolds
that may be defined in $\C^m$, or in Calabi--Yau $m$-folds, or more
generally in {\it almost Calabi--Yau $m$-folds}. They are calibrated
with respect to the $m$-form $\Re\Om$. They are fairly rigid and
well-behaved, so that compact SL $m$-folds $N$ occur in smooth moduli
spaces of dimension $b^1(N)$, for instance. They are important
in String Theory, and are expected to play a r\^ole in the
eventual explanation of Mirror Symmetry.

This article is intended as an introduction to special Lagrangian
geometry, and a survey of the author's research on the singularities
of SL $m$-folds, of directions in which the subject might develop in
the next few years, and of possible applications of it to Mirror
Symmetry and the SYZ Conjecture.

Sections \ref{bp2} and \ref{bp3} discuss general properties of
special Lagrangian submanifolds of $\C^m$, and ways to construct
examples. Then Section~\ref{bp4} defines Calabi--Yau and almost
Calabi--Yau manifolds, and their special Lagrangian submanifolds.
Section \ref{bp5} discusses the deformation and obstruction theory of
{\it compact}\/ SL $m$-folds, and properties of their moduli spaces.

In Section~\ref{bp6} we describe a theory of {\it isolated conical
singularities} in compact SL $m$-folds. Finally, Section~\ref{bp7}
briefly introduces String Theory, Mirror Symmetry and the {\it SYZ
Conjecture}, a conjectural explanation of Mirror Symmetry of
Calabi--Yau 3-folds, and discusses mathematical progress towards
clarifying and proving the conjecture.
\medskip

\noindent{\it Acknowledgments.} This article is based on lecture
courses given to the Summer School on Symplectic Geometry in
Nordfjordeid, Norway, in June 2001, and to the Clay Institute's
Summer School on the Global Theory of Minimal Surfaces, MSRI,
California, in July 2001. I would like to thank the organizers
of these for inviting me to speak.

Many people have helped me develop my ideas on special Lagrangian
geometry; amongst them I would particularly like to thank Robert 
Bryant, Mark Gross, Mark Haskins, Nigel Hitchin, Ian McIntosh,
Richard Thomas, and Karen Uhlenbeck.

\section{Special Lagrangian submanifolds in $\C^m$}
\label{bp2}

We begin by defining {\it calibrations} and {\it calibrated 
submanifolds}, following Harvey and Lawson~\cite{HaLa}.

\begin{dfn} Let $(M,g)$ be a Riemannian manifold. An {\it oriented
tangent $k$-plane} $V$ on $M$ is a vector subspace $V$ of
some tangent space $T_xM$ to $M$ with $\dim V=k$, equipped
with an orientation. If $V$ is an oriented tangent $k$-plane
on $M$ then $g\vert_V$ is a Euclidean metric on $V\!$, so 
combining $g\vert_V$ with the orientation on $V$ gives a 
natural {\it volume form} $\vol_V$ on $V\!$, which is a 
$k$-form on~$V\!$.

Now let $\vp$ be a closed $k$-form on $M$. We say that
$\vp$ is a {\it calibration} on $M$ if for every oriented
$k$-plane $V$ on $M$ we have $\vp\vert_V\le \vol_V\!$. Here
$\vp\vert_V=\al\cdot\vol_V$ for some $\al\in\R$, and 
$\vp\vert_V\le\vol_V$ if $\al\le 1$. Let $N$ be an 
oriented submanifold of $M$ with dimension $k$. Then 
each tangent space $T_xN$ for $x\in N$ is an oriented
tangent $k$-plane. We say that $N$ is a {\it calibrated 
\label{bp2def1}
submanifold} if $\vp\vert_{T_xN}=\vol_{T_xN}$ for all~$x\in N$.
\end{dfn}

It is easy to show that calibrated submanifolds are automatically
{\it minimal submanifolds} \cite[Th.~II.4.2]{HaLa}. Here is the 
definition of special Lagrangian submanifolds in $\C^m$, taken
from~\cite[\S III]{HaLa}.

\begin{dfn} Let $\C^m\cong\R^{2m}$ have complex coordinates 
$(z_1,\dots,z_m)$ and complex structure $I$, and define a metric 
$g'$, K\"ahler form $\om'$ and complex volume form $\Om'$ on $\C^m$ by
\e
\begin{split}
g'=\ms{\d z_1}+\cdots+\ms{\d z_m},\quad
\om'&=\frac{i}{2}(\d z_1\w\d\bar z_1+\cdots+\d z_m\w\d\bar z_m),\\
\text{and}\quad\Om'&=\d z_1\w\cdots\w\d z_m.
\end{split}
\label{bp2eq1}
\e
Then $\Re\Om'$ and $\Im\Om'$ are real $m$-forms on $\C^m$. Let
$L$ be an oriented real submanifold of $\C^m$ of real dimension 
$m$. We call $L$ a {\it special Lagrangian submanifold\/} in $\C^m$,
or {\it SL $m$-fold\/} for short, if $L$ is calibrated with respect
\label{bp2def2}
to $\Re\Om'$, in the sense of Definition~\ref{bp2def1}.
\end{dfn}

In fact there is a more general definition involving a {\it phase}
${\rm e}^{i\th}$: if $\th\in[0,2\pi)$, we say that $L$ is {\it special
Lagrangian with phase} ${\rm e}^{i\th}$ if it is calibrated with
respect to $\cos\th\,\Re\Om'+\sin\th\,\Im\Om'$. But we will not use
this.

We shall identify the family $\cal F$ of tangent $m$-planes in $\C^m$
calibrated with respect to $\Re\Om'$. The subgroup of $\GL(2m,\R)$ 
preserving $g',\om'$ and $\Om'$ is the Lie group $\SU(m)$ of complex 
unitary matrices with determinant 1. Define a real vector subspace 
$U$ in $\C^m$ to be
\e
U=\bigl\{(x_1,\ldots,x_m):x_j\in\R\bigr\}\subset\C^m,
\label{bp2eq2}
\e
and let $U$ have the usual orientation. Then $U$ is calibrated
with respect to $\Re\Om'$. 

Furthermore, any oriented real vector subspace $V$ in $\C^m$ calibrated 
with respect to  $\Re\Om'$ is of the form $V=\ga\cdot U$ for some $\ga\in\SU(m)$.
Therefore $\SU(m)$ acts transitively on $\cal F$. The stabilizer subgroup
of $U$ in $\SU(m)$ is the subset of matrices in $\SU(m)$ with real entries,
which is $\SO(m)$. Thus ${\cal F}\cong\SU(m)/\SO(m)$, and we prove:

\begin{prop} The family $\cal F$ of oriented real\/ $m$-dimensional
vector subspaces $V$ in $\C^m$ with\/ $\Re\Om'\vert_V=\vol_V$ is
\label{bp2prop1}
isomorphic to $\SU(m)/\SO(m)$, and has dimension~$\ha(m^2+m-2)$.
\end{prop}

The dimension follows because $\dim\SU(m)=m^2-1$ and $\dim\SO(m)=\ha m(m-1)$.
It is easy to see that $\om'\vert_U=\Im\Om'\vert_U=0$. As $\SU(m)$
preserves $\om'$ and $\Im\Om'$ and acts transitively on $\cal F$, it
follows that $\om'\vert_V=\Im\Om'\vert_V=0$ for any $V\in{\cal F}$.
Conversely, if $V$ is a real $m$-dimensional vector subspace of $\C^m$ 
and $\om'\vert_V=\Im\Om'\vert_V=0$, then $V$ lies in $\cal F$, with some 
orientation. This implies an alternative characterization of special 
Lagrangian submanifolds, \cite[Cor.~III.1.11]{HaLa}:

\begin{prop} Let\/ $L$ be a real\/ $m$-dimensional submanifold 
of\/ $\C^m$. Then $L$ admits an orientation making it into a
special Lagrangian submanifold of\/ $\C^m$ if and only if\/ 
\label{bp2prop2}
$\om'\vert_L\equiv 0$ and\/~$\Im\Om'\vert_L\equiv 0$.
\end{prop}

Note that an $m$-dimensional submanifold $L$ in $\C^m$ is 
called {\it Lagrangian} if $\om'\vert_L\equiv 0$. (This is a term
from symplectic geometry, and $\om'$ is a symplectic structure.) 
Thus special Lagrangian submanifolds are Lagrangian submanifolds 
satisfying the extra condition that $\Im\Om'\vert_L\equiv 0$, which 
is how they get their name.

\subsection{Special Lagrangian 2-folds in $\C^2$ and the quaternions}
\label{bp21}

The smallest interesting dimension, $m=2$, is a special case.
Let $\C^2$ have complex coordinates $(z_1,z_2)$, complex
structure $I$, and metric $g'$, K\"ahler form $\om'$ and holomorphic 
2-form $\Om'$ defined in \eq{bp2eq1}. Define real coordinates 
$(x_0,x_1,x_2,x_3)$ on $\C^2\cong\R^4$ by $z_0=x_0+ix_1$, 
$z_1=x_2+ix_3$. Then
\begin{alignat*}{2}
g'&=\d x_0^2+\cdots+\d x_3^2,&\qquad
\om'&=\d x_0\w\d x_1+\d x_2\w\d x_3,\\
\Re\Om'&=\d x_0\w\d x_2-\d x_1\w\d x_3&\quad\text{and}\quad
\Im\Om'&=\d x_0\w\d x_3+\d x_1\w\d x_2.
\end{alignat*}
Now define a {\it different} set of complex coordinates $(w_1,w_2)$ 
on $\C^2=\R^4$ by $w_1=x_0+ix_2$ and $w_2=x_1-ix_3$. 
Then~$\om'-i\Im\Om'=\d w_1\w\d w_2$. 

But by Proposition \ref{bp2prop2}, a real 2-submanifold $L\subset\R^4$ 
is special Lagrangian if and only if $\om'\vert_L\equiv\Im\Om'\vert_L
\equiv 0$. Thus, $L$ is special Lagrangian if and only if 
$(\d w_1\w\d w_2)\vert_L\equiv 0$. But this holds if and only if
$L$ is a {\it holomorphic curve} with respect to the complex 
coordinates~$(w_1,w_2)$. 

Here is another way to say this. There are {\it two different}
complex structures $I$ and $J$ involved in this problem,
associated to the two different complex coordinate systems 
$(z_1,z_2)$ and $(w_1,w_2)$ on $\R^4$. In the coordinates 
$(x_0,\ldots,x_3)$, $I$ and $J$ are given by
\begin{alignat*}{4}
I\Bigl({\frac{\pd}{\pd x_0}}\Bigr)&={\frac{\pd}{\pd x_1}},\quad &
I\Bigl({\frac{\pd}{\pd x_1}}\Bigr)&=-{\frac{\pd}{\pd x_0}},\quad &
I\Bigl({\frac{\pd}{\pd x_2}}\Bigr)&={\frac{\pd}{\pd x_3}},\quad &
I\Bigl({\frac{\pd}{\pd x_3}}\Bigr)&=-{\frac{\pd}{\pd x_2}},\\
J\Bigl({\frac{\pd}{\pd x_0}}\Bigr)&={\frac{\pd}{\pd x_2}},\quad &
J\Bigl({\frac{\pd}{\pd x_1}}\Bigr)&=-{\frac{\pd}{\pd x_3}},\quad &
J\Bigl({\frac{\pd}{\pd x_2}}\Bigr)&=-{\frac{\pd}{\pd x_0}},\quad &
J\Bigl({\frac{\pd}{\pd x_3}}\Bigr)&={\frac{\pd}{\pd x_1}}.
\end{alignat*}
The usual complex structure on $\C^2$ is $I$, but a 2-fold $L$ in
$\C^2$ is special Lagrangian if and only if it is holomorphic
with respect to  the alternative complex structure $J$. This means that 
special Lagrangian 2-folds are already very well understood,
so we generally focus our attention on dimensions~$m\ge 3$.

We can express all this in terms of the {\it quaternions} $\H$\,.
The complex structures $I,J$ anticommute, so that $IJ=-JI$, and 
$K=IJ$ is also a complex structure on $\R^4$, and $\an{1,I,J,K}$ 
is an algebra of automorphisms of $\R^4$ isomorphic to~$\H$\,.

\subsection{Special Lagrangian submanifolds in $\C^m$ as graphs}
\label{bp22}

In symplectic geometry, there is a well-known way of manufacturing 
{\it Lagrangian} submanifolds of $\R^{2m}\cong\C^m$, which works as 
follows. Let $f:\R^m\ra\R$ be a smooth function, and define
\begin{equation*}
\Ga_f\!=\!\bigl\{\bigl(x_1\!+\!i{\ts\frac{\pd f}{\pd x_1}}(x_1,\ldots,x_m),
\ldots,x_m\!+\!i{\ts\frac{\pd f}{\pd x_m}}(x_1,\ldots,x_m)\bigr):
x_1,\ldots,x_m\!\in\!\R\bigr\}.
\end{equation*}
Then $\Ga_f$ is a smooth real $m$-dimensional submanifold of $\C^m$,
with $\om'\vert_{\Ga_f}\equiv 0$. Identifying $\C^m\cong\R^{2m}\cong
\R^m\t(\R^m)^*$, we may regard $\Ga_f$ as the graph of the 1-form
$\d f$ on $\R^m$, so that $\Ga_f$ is the {\it graph of a closed\/ 
$1$-form}. Locally, but not globally, every Lagrangian submanifold 
arises from this construction.

Now by Proposition \ref{bp2prop2}, a special Lagrangian $m$-fold
in $\C^m$ is a Lagrangian $m$-fold $L$ satisfying the additional
condition that $\Im\Om'\vert_L\equiv 0$. We shall find the
condition for $\Ga_f$ to be a special Lagrangian $m$-fold. Define
the {\it Hessian} $\Hess f$ of $f$ to be the $m\t m$ matrix 
$\bigl(\frac{\pd^2f}{\pd x_i\pd x_j}\bigr)_{i,j=1}^m$ of real 
functions on $\R^m$. Then it is easy to show that $\Im\Om'
\vert_{\Ga_f}\equiv 0$ if and only if
\e
\ts\Im\det_{\sst\mathbb C}\bigl(I+i\Hess f\bigr)\equiv 0
\quad\text{on $\C^m$.}
\label{bp2eq3}
\e
This is a {\it nonlinear second-order elliptic partial differential
equation} upon the function~$f:\R^m\ra\R$\,. 

\subsection{Local discussion of special Lagrangian deformations}
\label{bp23}

Suppose $L_0$ is a special Lagrangian submanifold in $\C^m$ (or, more
generally, in some (almost) Calabi--Yau $m$-fold). What can we say
about the family of {\it special Lagrangian deformations} of $L_0$,
that is, the set of special Lagrangian $m$-folds $L$ that are
``close to $L_0$'' in a suitable sense? Essentially, deformation
theory is one way of thinking about the question ``how many special 
Lagrangian submanifolds are there in~$\C^m$''?

Locally (that is, in small enough open sets), every special Lagrangian
$m$-fold looks quite like $\R^m$ in $\C^m$. Therefore deformations of 
special Lagrangian $m$-folds should look like special Lagrangian 
deformations of $\R^m$ in $\C^m$. So, we would like to know what
special Lagrangian $m$-folds $L$ in $\C^m$ close to $\R^m$ look like.

Now $\R^m$ is the graph $\Ga_f$ associated to the function $f\equiv 0$.
Thus, a graph $\Ga_f$ will be close to $\R^m$ if the function $f$ and 
its derivatives are small. But then $\Hess f$ is small, so we can 
approximate equation \eq{bp2eq3} by its {\it linearization}. For
\begin{equation*}
\ts\Im\det_{\sst\mathbb C}\bigl(I+i\Hess f\bigr)=
\Tr\Hess f+\text{higher order terms}.
\end{equation*}
Thus, when the second derivatives of $f$ are small, equation
\eq{bp2eq3} reduces approximately to $\Tr\Hess f\equiv 0$. But 
\begin{equation*}
Tr\Hess f=\frac{\pd^2f}{(\pd x_1)^2}+\cdots+\frac{\pd^2f}{(\pd x_m)^2}
=-\De f,
\end{equation*}
where $\De$ is the {\it Laplacian} on~$\R^m$.

Hence, the small special Lagrangian deformations of $\R^m$ in
$\C^m$ are approximately parametrized by small {\it harmonic 
functions} on $\R^m$. Actually, because adding a constant to
$f$ has no effect on $\Ga_f$, this parametrization is degenerate.
We can get round this by parametrizing instead by $\d f$, which
is a closed and coclosed 1-form. This justifies the following:
\medskip

\noindent{\bf Principle.} {\it Small special Lagrangian deformations 
of a special Lagrangian $m$-fold\/ $L$ are approximately parametrized 
by closed and coclosed\/ $1$-forms $\al$ on~$L$.}
\medskip

\noindent This is the idea behind McLean's Theorem, Theorem 
\ref{bp5thm1} below.

We have seen using \eq{bp2eq3} that the deformation problem for special 
Lagrangian $m$-folds can be written as an {\it elliptic equation}.
In particular, there are the same number of equations as functions, so
the problem is neither overdetermined nor underdetermined. Therefore we
do not expect special Lagrangian $m$-folds to be very few and very rigid 
(as would be the case if \eq{bp2eq3} were overdetermined), nor to 
be very abundant and very flabby (as would be the case if 
\eq{bp2eq3} were underdetermined).

If we think about Proposition \ref{bp2prop1} for a while, this may
seem surprising. For the set $\cal F$ of special Lagrangian $m$-planes
in $\C^m$ has dimension $\ha(m^2+m-2)$, but the set of all real $m$-planes
in $\C^m$ has dimension $m^2$. So the special Lagrangian $m$-planes
have codimension $\ha(m^2-m+2)$ in the set of all $m$-planes.

This means that the condition for a real $m$-submanifold $L$ in
$\C^m$ to be special Lagrangian is $\ha(m^2-m+2)$ real equations on 
each tangent space of $L$. However, the freedom to vary $L$ is the
sections of its normal bundle in $\C^m$, which is $m$ real functions. 
When $m\ge 3$, there are more equations than functions, so we would 
expect the deformation problem to be {\it overdetermined}.

The explanation is that because $\om'$ is a {\it closed}\/ 2-form,
submanifolds $L$ with $\om'\vert_L\equiv 0$ are much more abundant
than would otherwise be the case. So the closure of $\om'$ is a
kind of integrability condition necessary for the existence of 
many special Lagrangian submanifolds, just as the integrability 
of an almost complex structure is a necessary condition for the
existence of many complex submanifolds of dimension greater than 
1 in a complex manifold.

\section{Constructions of SL $m$-folds in $\C^m$}
\label{bp3}

We now describe five methods of constructing special Lagrangian
$m$-folds in $\C^m$, drawn from papers by the author
\cite{Joyc2,Joyc3,Joyc5,Joyc6,Joyc7,Joyc8,Joyc9,Joyc10},
Bryant \cite{Brya}, Castro and Urbano \cite{CaUr}, Goldstein
\cite{Gold1,Gold2}, Harvey \cite[p.~139--143]{Harv}, Harvey and
Lawson \cite[\S III]{HaLa}, Haskins \cite{Hask}, Lawlor
\cite{Lawl}, Ma and Ma \cite{MaMa}, McIntosh \cite{McIn} and
Sharipov \cite{Shar}. These yield many examples of singular SL
$m$-folds, and so hopefully will help in understanding
general singularities of SL $m$-folds in Calabi--Yau $m$-folds.

\subsection{SL $m$-folds with large symmetry groups}
\label{bp31}

Here is a method used in \cite{Joyc5} (and also by Harvey and
Lawson \cite[\S III.3]{HaLa}, Haskins \cite{Hask} and Goldstein
\cite{Gold1,Gold2}) to construct examples of SL $m$-folds in 
$\C^m$. The group $\SU(m)\lt\C^m$ acts on $\C^m$ preserving 
all the structure $g',\om',\Om'$, so that it takes SL $m$-folds 
to SL $m$-folds in $\C^m$. Let $G$ be a Lie subgroup of 
$\SU(m)\lt\C^m$ with Lie algebra $\g$, and $N$ a connected 
$G$-invariant SL $m$-fold in~$\C^m$.

Since $G$ preserves the symplectic form $\om'$ on $\C^m$, one 
can show that it has a {\it moment map} $\mu:\C^m\ra\g^*$.
As $N$ is Lagrangian, one can show that $\mu$ is constant
on $N$, that is, $\mu\equiv c$ on $N$ for some $c\in Z(\g^*)$, 
the {\it center} of~$\g^*$.

If the orbits of $G$ in $N$ are of codimension 1 (that is,
dimension $m-1$), then $N$ is a 1-parameter family of 
$G$-orbits ${\cal O}_t$ for $t\in\R$\,. After reparametrizing 
the variable $t$, it can be shown that the special Lagrangian 
condition is equivalent to an ODE in $t$ upon the 
orbits~${\cal O}_t$.

Thus, we can construct examples of cohomogeneity one SL $m$-folds
in $\C^m$ by solving an ODE in the family of $(m-1)$-dimensional 
$G$-orbits $\cal O$ in $\C^m$ with $\mu\vert_{\cal O}\equiv c$,
for fixed $c\in Z(\g^*)$. This ODE usually turns out to be 
{\it integrable}.

Now suppose $N$ is a {\it special Lagrangian cone} in $\C^m$, 
invariant under a subgroup $G\subset\SU(m)$ which has orbits of 
dimension $m-2$ in $N$. In effect the symmetry group of $N$ is 
$G\t\R_+$, where $\R_+$ acts by {\it dilations}, as $N$ is a cone. 
Thus, in this situation too the symmetry group of $N$ acts with 
cohomogeneity one, and we again expect the problem to reduce to 
an ODE.

One can show that $N\cap{\cal S}^{2m-1}$ is a 1-parameter family of 
$G$-orbits ${\cal O}_t$ in ${\cal S}^{2m-1}\cap\mu^{-1}(0)$ satisfying 
an ODE. By solving this ODE we construct SL cones in $\C^m$.
When $G=\U(1)^{m-2}$, the ODE has many {\it periodic solutions} 
which give large families of distinct SL cones on $T^{m-1}$. In 
particular, we can find many examples of SL $T^2$-cones in~$\C^3$.

\subsection{Evolution equations for SL $m$-folds}
\label{bp32}

The following method was used in \cite{Joyc2} and \cite{Joyc3} to 
construct many examples of SL $m$-folds in $\C^m$. A related but 
less general method was used by Lawlor \cite{Lawl}, and completed 
by Harvey~\cite[p.~139--143]{Harv}.

Let $P$ be a real analytic $(m-1)$-dimensional manifold, and 
$\chi$ a nonvanishing real analytic section of $\La^{m-1}TP$.
Let $\{\phi_t:t\in\R\}$ be a 1-parameter family of real analytic 
maps $\phi_t:P\ra\C^m$. Consider the ODE
\e
\Bigl(\frac{\d\phi_t}{\d t}\Bigr)^{\!b}=(\phi_t)_*(\chi)^{a_1\ldots 
a_{m-1}}(\Re\Om')_{a_1\ldots a_{m-1}a_m}g'{}^{a_mb},
\label{bp3eq1}
\e
using the index notation for (real) tensors on $\C^m$, where 
$g'{}^{ab}$ is the inverse of the Euclidean metric $g'_{ab}$
on~$\C^m$.

It is shown in \cite[\S 3]{Joyc2} that if the $\phi_t$ satisfy 
\eq{bp3eq1} and $\phi_0^*(\om')\equiv 0$, then $\phi_t^*(\om')\equiv 0$ 
for all $t$, and $N=\bigl\{\phi_t(p):p\in P$, $t\in\R\bigr\}$ is 
an SL $m$-fold in $\C^m$ wherever it is nonsingular. We think of 
\eq{bp3eq1} as an {\it evolution equation}, and $N$ as the result 
of evolving a 1-parameter family of $(m\!-\!1)$-submanifolds
$\phi_t(P)$ in~$\C^m$.

Here is one way to understand this result. Suppose we are given 
$\phi_t:P\ra\C^m$ for some $t$, and we want to find an SL $m$-fold
$N$ in $\C^m$ containing the $(m\!-\!1)$-submanifold $\phi_t(P)$. 
As $N$ is Lagrangian, a necessary condition for this is that $\om'
\vert_{\phi_t(P)}\equiv 0$, and hence $\phi_t^*(\om')\equiv 0$ on~$P$. 

The effect of equation \eq{bp3eq1} is to flow $\phi_t(P)$ in the 
direction in which $\Re\Om'$ is ``largest''. The result is that $\Re\Om'$ 
is ``maximized'' on $N$, given the initial conditions. But $\Re\Om'$ is
maximal on $N$ exactly when $N$ is calibrated with respect to $\Re\Om'$, that
is, when $N$ is special Lagrangian. The same technique also works for 
other calibrations, such as the associative and coassociative 
calibrations on $\R^7$, and the Cayley calibration on~$\R^8$.

Now \eq{bp3eq1} evolves amongst the infinite-dimensional family
of real analytic maps $\phi:P\ra\C^m$ with $\phi^*(\om')\equiv 0$,
so it is an {\it infinite-dimensional} problem, and thus difficult
to solve explicitly. However, there are {\it finite-dimensional} 
families $\cal C$ of maps $\phi:P\ra\C^m$ such that evolution 
stays in $\cal C$. This gives a {\it finite-dimensional} ODE, 
which can hopefully be solved fairly explicitly. For example, 
if we take $G$ to be a Lie subgroup of $\SU(m)\lt\C^m$, $P$ to 
be an $(m\!-\!1)$-dimensional homogeneous space $G/H$, and 
$\phi:P\ra\C^m$ to be $G$-equivariant, we recover the 
construction of~Section~\ref{bp31}. 

But there are also other possibilities for $\cal C$ which do not 
involve a symmetry assumption. Suppose $P$ is a submanifold of 
$\R^n$, and $\chi$ the restriction to $P$ of a linear or affine 
map $\R^n\ra\La^{m-1}\R^n$. (This is a strong condition on $P$
and $\chi$.) Then we can take $\cal C$ to be the set of 
restrictions to $P$ of linear or affine maps~$\R^n\ra\C^m$.

For instance, set $m=n$ and let $P$ be a quadric in $\R^m$. Then
one can construct SL $m$-folds in $\C^m$ with few symmetries by
evolving quadrics in Lagrangian planes $\R^m$ in $\C^m$. When $P$ 
is a quadric cone in $\R^m$ this gives many SL cones on products
of spheres~${\cal S}^a\t{\cal S}^b\t{\cal S}^1$.

\subsection{Ruled special Lagrangian 3-folds}
\label{bp33}

A 3-submanifold $N$ in $\C^3$ is called {\it ruled} if it is fibered by 
a 2-dimensional family $\cal F$ of real lines in $\C^3$. A {\it cone} 
$N_0$ in $\C^3$ is called {\it two-sided} if $N_0=-N_0$. Two-sided 
cones are automatically ruled. If $N$ is a ruled 3-fold in $\C^3$, 
we define the {\it asymptotic cone} $N_0$ of $N$ to be the two-sided 
cone fibered by the lines passing through 0 and parallel to those 
in~$\cal F$. 

Ruled SL 3-folds are studied in \cite{Joyc6}, and also by
Harvey and Lawson \cite[\S III.3.C, \S III.4.B]{HaLa} and
Bryant \cite[\S 3]{Brya}. Each (oriented) real line in $\C^3$ 
is determined by its {\it direction} in ${\cal S}^5$ together 
with an orthogonal {\it translation} from the origin. Thus a 
ruled 3-fold $N$ is determined by a 2-dimensional family of 
directions and translations. 

The condition for $N$ to be special Lagrangian turns out 
\cite[\S 5]{Joyc6} to reduce to two equations, the first 
involving only the direction components, and the second 
{\it linear} in the translation components. Hence, if a 
ruled 3-fold $N$ in $\C^3$ is special Lagrangian, then so is 
its asymptotic cone $N_0$. Conversely, the ruled SL 3-folds 
$N$ asymptotic to a given two-sided SL cone $N_0$ come from 
solutions of a linear equation, and so form a {\it vector space}.

Let $N_0$ be a two-sided SL cone, and let $\Si=N_0\cap{\cal S}^5$. 
Then $\Si$ is a {\it Riemann surface}. Holomorphic vector fields on 
$\Si$ give solutions to the linear equation (though not all solutions)
\cite[\S 6]{Joyc6}, and so yield new ruled SL 3-folds. In particular, 
each SL $T^2$-cone gives a 2-dimensional family of ruled SL 3-folds, 
which are generically diffeomorphic to $T^2\t\R$ as immersed 
3-submanifolds. 

\subsection{Integrable systems}
\label{bp34}

Let $N_0$ be a special Lagrangian cone in $\C^3$, and set
$\Si=N_0\cap{\cal S}^5$. As $N_0$ is calibrated, it is minimal
in $\C^3$, and so $\Si$ is minimal in ${\cal S}^5$. That is, 
$\Si$ is a {\it minimal Legendrian surface} in ${\cal S}^5$.
Let $\pi:{\cal S}^5\ra\CP^2$ be the Hopf projection. One can 
also show that $\pi(\Si)$ is a {\it minimal Lagrangian surface}
in~$\CP^2$.

Regard $\Si$ as a {\it Riemann surface}. Then the inclusions 
$\iota:\Si\ra{\cal S}^5$ and $\pi\circ\iota:\Si\ra\CP^2$ 
are {\it conformal harmonic maps}. Now harmonic maps from 
Riemann surfaces into ${\cal S}^n$ and $\CP^m$ are an 
{\it integrable system}. There is a complicated theory 
for classifying them in terms of algebro-geometric ``spectral 
data'', and finding ``explicit'' solutions. In principle, this
gives all harmonic maps from $T^2$ into ${\cal S}^n$ and
$\CP^m$. So, the field of integrable systems offers the hope
of a {\it classification} of all SL $T^2$-cones in~$\C^3$.

For a good general introduction to this field, see Fordy and
Wood \cite{FoWo}. Sharipov \cite{Shar} and Ma and Ma \cite{MaMa}
apply this integrable systems machinery to describe minimal
Legendrian tori in ${\cal S}^5$, and minimal Lagrangian tori
in $\CP^2$, respectively, giving explicit formulae in terms
of Prym theta functions. McIntosh \cite{McIn} provides a more
recent, readable, and complete discussion of special Lagrangian
cones in $\C^3$ from the integrable systems perspective.

The families of SL $T^2$-cones constructed by $\U(1)$-invariance 
in Section~\ref{bp31}, and by evolving quadrics in Section~\ref{bp32}, turn out
to come from a more general, very explicit, ``integrable systems'' family 
of conformal harmonic maps $\R^2\ra{\cal S}^5$ with Legendrian image, 
involving two commuting, integrable ODEs, described in \cite{Joyc7}.
So, we can fit some of our examples into the integrable systems framework.

However, we know a good number of other constructions of SL $m$-folds 
in $\C^m$ which have the classic hallmarks of integrable systems
--- elliptic functions, commuting ODEs, and so on --- but which
are not yet understood from the point of view of integrable systems.

\subsection{Analysis and $\U(1)$-invariant SL $3$-folds in $\C^3$}
\label{bp35}

Next we summarize the author's three papers \cite{Joyc8,Joyc9,Joyc10},
which study SL 3-folds $N$ in $\C^3$ invariant under the $\U(1)$-action
\e
{\rm e}^{i\th}:(z_1,z_2,z_3)\mapsto
({\rm e}^{i\th}z_1,{\rm e}^{-i\th}z_2,z_3)
\quad\text{for ${\rm e}^{i\th}\in\U(1)$.}
\label{bp3eq2}
\e
These three papers are surveyed in \cite{Joyc11}. Locally we can
write $N$ in the form
\e
\begin{split}
N=\bigl\{(z_1,z_2,z_3)\in\C^3:\,& z_1z_2=v(x,y)+iy,\quad z_3=x+iu(x,y),\\
&\ms{z_1}-\ms{z_2}=2a,\quad (x,y)\in S\bigr\},
\end{split}
\label{bp3eq3}
\e
where $S$ is a domain in $\R^2$, $a\in\R$ and $u,v:S\ra\R$ are
continuous.

Here we may take $\ms{z_1}-\ms{z_2}=2a$ to be one of the equations
defining $N$ as $\ms{z_1}-\ms{z_2}$ is the {\it moment map} of the
$\U(1)$-action \eq{bp3eq2}, and so $\ms{z_1}-\ms{z_2}$ is constant
on any $\U(1)$-invariant Lagrangian 3-fold in $\C^3$. Effectively
\eq{bp3eq3} just means that we are choosing $x=\Re(z_3)$ and
$y=\Im(z_1z_2)$ as local coordinates on the 2-manifold $N/\U(1)$.
Then we find~\cite[Prop.~4.1]{Joyc8}:

\begin{prop} Let\/ $S,a,u,v$ and\/ $N$ be as above. Then
\label{bp3prop1}
\begin{itemize}
\item[{\rm(a)}] If\/ $a=0$, then $N$ is a (possibly singular) SL\/
$3$-fold in $\C^3$ if\/ $u,v$ are differentiable and satisfy
\e
\frac{\pd u}{\pd x}=\frac{\pd v}{\pd y}
\quad\text{and}\quad
\frac{\pd v}{\pd x}=-2\bigl(v^2+y^2\bigr)^{1/2}\frac{\pd u}{\pd y},
\label{bp3eq4}
\e
except at points $(x,0)$ in $S$ with\/ $v(x,0)=0$, where $u,v$ 
need not be differentiable. The singular points of\/ $N$ are those
of the form $(0,0,z_3)$, where $z_3=x+iu(x,0)$ for $(x,0)\in S$ 
with\/~$v(x,0)=0$.
\item[{\rm(b)}] If\/ $a\ne 0$, then $N$ is a nonsingular SL\/ $3$-fold
in $\C^3$ if and only if\/ $u,v$ are differentiable in $S$ and
satisfy
\e
\frac{\pd u}{\pd x}=\frac{\pd v}{\pd y}\quad\text{and}\quad
\frac{\pd v}{\pd x}=-2\bigl(v^2+y^2+a^2\bigr)^{1/2}\frac{\pd u}{\pd y}.
\label{bp3eq5}
\e
\end{itemize}
\end{prop}

Now \eq{bp3eq4} and \eq{bp3eq5} are {\it nonlinear Cauchy--Riemann
equations}. Thus, we may treat $u+iv$ as like a holomorphic function
of $x+iy$. Many of the results in \cite{Joyc8,Joyc9,Joyc10} are
analogues of well-known results in elementary complex analysis.

In \cite[Prop.~7.1]{Joyc8} we show that solutions $u,v\in C^1(S)$ 
of \eq{bp3eq5} come from a potential $f\in C^2(S)$ satisfying a
second-order quasilinear elliptic equation.

\begin{prop} Let\/ $S$ be a domain in $\R^2$ and\/ $u,v\in C^1(S)$
satisfy \eq{bp3eq5} for $a\ne 0$. Then there exists $f\in C^2(S)$
with\/ $\frac{\pd f}{\pd y}=u$, $\frac{\pd f}{\pd x}=v$ and
\e
P(f)=\Bigl(\Bigl(\frac{\pd f}{\pd x}\Bigr)^2+y^2+a^2
\Bigr)^{-1/2}\frac{\pd^2f}{\pd x^2}+2\,\frac{\pd^2f}{\pd y^2}=0.
\label{bp3eq6}
\e
This $f$ is unique up to addition of a constant, $f\mapsto f+c$.
Conversely, all solutions of\/ \eq{bp3eq6} yield solutions 
\label{bp3prop2}
of\/~\eq{bp3eq5}. 
\end{prop}

In the following result, a condensation of \cite[Th.~7.6]{Joyc8}
and \cite[Th.s 9.20 \& 9.21]{Joyc9}, we prove existence and
uniqueness for the {\it Dirichlet problem} for~\eq{bp3eq6}.

\begin{thm} Suppose $S$ is a strictly convex domain in $\R^2$ invariant
under $(x,y)\mapsto(x,-y)$, and\/ $\al\in(0,1)$. Let\/ $a\in\R$ and\/
$\phi\in C^{3,\al}(\pd S)$. Then if\/ $a\ne 0$ there exists a unique
solution $f$ of\/ \eq{bp3eq6} in $C^{3,\al}(S)$ with\/ $f\vert_{\pd S}
=\phi$. If\/ $a=0$ there exists a unique $f\in C^1(S)$ with\/
$f\vert_{\pd S}=\phi$, which is twice weakly differentiable and
satisfies \eq{bp3eq6} with weak derivatives. Furthermore, the map
$C^{3,\al}(\pd S)\t\R\ra C^1(S)$ taking $(\phi,a)\mapsto f$ is
\label{bp3thm1}
continuous.
\end{thm}

Here a domain $S$ in $\R^2$ is {\it strictly convex} if it is
convex and the curvature of $\pd S$ is nonzero at each point.
Also domains are by definition compact, with smooth boundary,
and $C^{3,\al}(\pd S)$ and $C^{3,\al}(S)$ are {\it H\"older
spaces} of functions on $\pd S$ and $S$. For more details
see~\cite{Joyc8,Joyc9}.

Combining Propositions \ref{bp3prop1} and \ref{bp3prop2} and
Theorem \ref{bp3thm1} gives existence and uniqueness for a large
class of $\U(1)$-invariant SL 3-folds in $\C^3$, with boundary
conditions, and including {\it singular} SL 3-folds. It is
interesting that this existence and uniqueness is {\it entirely
unaffected} by singularities appearing in~$S^\circ$. 

Here are some other areas covered in \cite{Joyc8,Joyc9,Joyc10}.
Examples of solutions $u,v$ of \eq{bp3eq4} and \eq{bp3eq5} are
given in \cite[\S 5]{Joyc8}. In \cite{Joyc9} we give more
precise statements on the regularity of singular solutions of
\eq{bp3eq4} and \eq{bp3eq6}. In \cite[\S 6]{Joyc8} and \cite[\S
7]{Joyc10} we consider the zeroes of $(u_1,v_1)-(u_2,v_2)$,
where $(u_j,v_j)$ are (possibly singular) solutions of \eq{bp3eq4}
and~\eq{bp3eq5}.

We show that if $(u_1,v_1)\not\equiv(u_2,v_2)$ then the zeroes
of $(u_1,v_1)-(u_2,v_2)$ in $S^\circ$ are {\it isolated}, with
a positive integer {\it multiplicity}, and that the zeroes of
$(u_1,v_1)-(u_2,v_2)$ in $S^\circ$ can be counted with
multiplicity in terms of boundary data on $\pd S$. In particular,
under some boundary conditions we can show $(u_1,v_1)-(u_2,v_2)$
has no zeroes in $S^\circ$, so that the corresponding SL 3-folds
do not intersect. This will be important in constructing
$\U(1)$-invariant SL fibrations in~Section~\ref{bp75}.

In \cite[\S 9--\S 10]{Joyc10} we study singularities of solutions
$u,v$ of \eq{bp3eq4}. We show that either $u(x,-y)\equiv u(x,y)$
and $v(x,-y)\equiv -v(x,y)$, so that $u,v$ are singular all
along the $x$-axis, or else the singular points of $u,v$ in
$S^\circ$ are all {\it isolated}, with a positive integer
{\it multiplicity}, and one of two {\it types}. We also show
that singularities exist with every multiplicity and type, and
multiplicity $n$ singularities occur in codimension $n$ in the
family of all $\U(1)$-invariant SL 3-folds.

\subsection{Examples of singular special Lagrangian 3-folds in $\C^3$}
\label{bp36}

We shall now describe four families of SL 3-folds in $\C^3$, as examples 
of the material of Sections~\ref{bp31}--\ref{bp34}. They have been chosen to 
illustrate different kinds of singular behavior of SL 3-folds, and 
also to show how nonsingular SL 3-folds can converge to a singular SL 
3-fold, to serve as a preparation for our discussion of singularities 
of SL $m$-folds in~Section~\ref{bp6}.

Our first example derives from Harvey and Lawson \cite[\S III.3.A]{HaLa},
and is discussed in detail in~\cite[\S 3]{Joyc4}.

\begin{ex} Define a subset $L_0$ in $\C^3$ by
\begin{equation*}
L_0=\bigl\{(r{\rm e}^{i\th_1},r{\rm e}^{i\th_2},r{\rm e}^{i\th_3}):
r\ge 0,\quad \th_1,\th_2,\th_3\in\R,\quad \th_1+\th_2+\th_3=0\bigr\}.
\end{equation*}
Then $L_0$ is a {\it special Lagrangian cone} on $T^2$. An 
alternative definition is
\begin{equation*}
L_0=\bigl\{(z_1,z_2,z_3)\in\C^3:\md{z_1}=\md{z_2}=\md{z_3},\;
\Im(z_1z_2z_3)=0,\; \Re(z_1z_2z_3)\ge 0\bigr\}.
\end{equation*}

Let $t>0$, write ${\cal S}^1=\bigl\{{\rm e}^{i\th}:\th\in\R\bigr\}$, 
and define a map $\phi_t:{\cal S}^1\t\C\ra\C^3$ by
\begin{equation*}
\phi_t:(e^{i\th},z)\mapsto
\bigl((\md{z}^2+t^2)^{1/2}{\rm e}^{i\th},z,e^{-i\th}\bar z\bigr).
\end{equation*}
Then $\phi_t$ is an {\it embedding}. Define $L_t=\Image\phi_t$.
Then $L_t$ is a nonsingular special Lagrangian 3-fold in $\C^3$
diffeomorphic to ${\cal S}^1\t\R^2$. An equivalent definition~is
\begin{align*}
L_t=\bigl\{(z_1,z_2,z_3)\in\C^3:\,&\ms{z_1}-t^2=\ms{z_2}=\ms{z_3},\\
&\Im(z_1z_2z_3)=0,\quad \Re(z_1z_2z_3)\ge 0\bigr\}.
\end{align*}

As $t\ra 0_+$, the nonsingular SL 3-fold $L_t$ converges to the
singular SL cone $L_0$. Note that $L_t$ is {\it asymptotic}
to $L_0$ at infinity, and that $L_t=t\,L_1$ for $t>0$, so that
the $L_t$ for $t>0$ are all homothetic to each other. Also,
each $L_t$ for $t\ge 0$ is invariant under the $T^2$ subgroup 
of $\SU(3)$ acting by
\begin{equation*}
(z_1,z_2,z_3)\mapsto({\rm e}^{i\th_1}z_1,{\rm e}^{i\th_2}z_2,
{\rm e}^{i\th_3}z_3) \
\text{for $\th_1,\th_2,\th_3\in\R$ with $\th_1+\th_2+\th_3=0$,}
\end{equation*}
and so fits into the framework of Section~\ref{bp31}. By
\cite[Th.~5.1]{Joyc8} the $L_a$ may also be written in the form
\label{bp3ex1}
\eq{bp3eq3} for continuous $u,v:\R^2\ra\R$, as in~Section~\ref{bp35}.
\end{ex}

Our second example is adapted from Harvey and 
Lawson~\cite[\S III.3.B]{HaLa}. 

\begin{ex} For each $t>0$, define
\begin{equation*}
L_t=\bigl\{({\rm e}^{i\th}x_1,{\rm e}^{i\th}x_2,{\rm e}^{i\th}x_3):
x_j\in\R,\ \th\in(0,\pi/3),\
x_1^2\!+\!x_2^2\!+\!x_3^2\!=\!t^2(\sin 3\th)^{-2/3}\bigr\}.
\end{equation*}
Then $L_t$ is a nonsingular embedded SL 3-fold in $\C^3$ diffeomorphic
to ${\cal S}^2\t\R$\,. As $t\ra 0_+$ it converges to the singular union
$L_0$ of the two SL 3-planes
\begin{equation*}
\Pi_1=\bigl\{(x_1,x_2,x_3):x_j\in\R\bigr\}\;\>\text{and}\;\>
\Pi_2=\bigl\{({\rm e}^{i\pi/3}x_1,{\rm e}^{i\pi/3}x_2,
{\rm e}^{i\pi/3}x_3):x_j\in\R\bigr\},
\end{equation*}
which intersect at 0. Note that $L_t$ is invariant under the action 
of the Lie subgroup $\SO(3)$ of $\SU(3)$, acting on $\C^3$ in the 
obvious way, so again this comes from the method of Section~\ref{bp31}. 
\label{bp3ex2}
Also $L_t$ is asymptotic to $L_0$ at infinity.
\end{ex}

Our third example is taken from~\cite[Ex.~9.4 \& Ex.~9.5]{Joyc5}.

\begin{ex} Let $a_1,a_2$ be positive, coprime integers, and set
$a_3=-a_1-a_2$. Let $c\in\R$, and define
\begin{equation*}
L^{a_1,a_2}_c=\bigl\{({\rm e}^{ia_1\th}x_1,{\rm e}^{ia_2\th}x_2,
i{\rm e}^{ia_3\th}x_3):\th\in\R,\; x_j\in\R,\; 
a_1x_1^2+a_2x_2^2+a_3x_3^2=c\bigr\}.
\end{equation*}
Then $L^{a_1,a_2}_c$ is an SL 3-fold, which comes from the
``evolving quadrics'' construction of Section~\ref{bp32}. It
is also symmetric under the $\U(1)$-action
\begin{equation*}
(z_1,z_2,z_3)\mapsto({\rm e}^{ia_1\th}z_1,{\rm e}^{ia_2\th}z_2,
i{\rm e}^{ia_3\th}z_3)\quad\text{for $\th\in\R$,}
\end{equation*}
but this is not a necessary feature of the construction; these
are just the easiest examples to write down.

When $c=0$ and $a_3$ is odd, $L^{a_1,a_2}_0$ is an embedded 
special Lagrangian cone on $T^2$, with one singular point at 0. 
When $c=0$ and $a_3$ is even, $L^{a_1,a_2}_0$ is two opposite 
embedded SL $T^2$-cones with one singular point at~0.

When $c>0$ and $a_3$ is odd, $L^{a_1,a_2}_c$ is an embedded 3-fold 
diffeomorphic to a nontrivial real line bundle over the Klein bottle.
When $c>0$ and $a_3$ is even, $L^{a_1,a_2}_c$ is an embedded 3-fold 
diffeomorphic to $T^2\t\R$\,. In both cases, $L^{a_1,a_2}_c$ is
a {\it ruled}\/ SL 3-fold, as in Section~\ref{bp33}, since it is fibered
by hyperboloids of one sheet in $\R^3$, which are ruled in two
different ways.

When $c<0$ and $a_3$ is odd, $L^{a_1,a_2}_c$ an immersed
copy of ${\cal S}^1\t\R^2$. When $c<0$ and $a_3$ is even, 
\label{bp3ex3}
$L^{a_1,a_2}_c$ two immersed copies of~${\cal S}^1\t\R^2$. 
\end{ex}

All the singular SL 3-folds we have seen so far have been {\it 
cones} in $\C^3$. Our final example, taken from \cite{Joyc3},
has more complicated singularities which are not cones. They
are difficult to describe in a simple way, so we will not say
much about them. For more details, see~\cite{Joyc3}.
                                                   
\begin{ex} In \cite[\S 5]{Joyc3} the author constructed a
family of maps $\Phi:\R^3\ra\C^3$ with special Lagrangian image
$N=\Image\Phi$. It is shown in \cite[\S 6]{Joyc3} that generic 
$\Phi$ in this family are immersions, so that $N$ is nonsingular
as an immersed SL 3-fold, but in codimension 1 in the family they 
develop isolated singularities.

Here is a rough description of these singularities, taken 
from \cite[\S 6]{Joyc3}. Taking the singular point to be at 
$\Phi(0,0,0)=0$, one can write $\Phi$ as
\e
\begin{split}
\Phi(x,y,t)=
&\bigl(x+{\ts\frac{1}{4}}g'({\bf u},{\bf v})t^2\bigr)\,{\bf u}
+\bigl(y^2-{\ts\frac{1}{4}}\ms{{\bf u}}t^2\bigr)\,{\bf v}\\
&+2yt\,{\bf u}\t{\bf v}+O\bigl(x^2+\md{xy}+\md{xt}+\md{y}^3+\md{t}^3\bigr),
\end{split}
\label{bp3eq7}
\e
where ${\bf u},{\bf v}$ are linearly independent vectors in $\C^3$ with 
$\om'({\bf u},{\bf v})=0$, and $\t:\C^3\t\C^3\ra\C^3$ is defined by
\begin{equation*}
(r_1,r_2,r_3)\t(s_1,s_2,s_3)={\ts\frac{1}{2}}(\bar r_2\bar s_3-\bar r_3\bar 
s_2,\bar r_3\bar s_1-\bar r_1\bar s_3,\bar r_1\bar s_2-\bar r_2\bar s_1).
\end{equation*}
The next few terms in the expansion \eq{bp3eq7} can also be given very 
explicitly, but we will not write them down as they are rather complex,
and involve further choices of vectors~${\bf w},{\bf x},\dots{}$. 

What is going on here is that the lowest order terms in $\Phi$ are a
{\it double cover} of the special Lagrangian plane $\an{{\bf u},{\bf v},
{\bf u}\t{\bf v}}_{\sst\mathbb R}$ in $\C^3$, {\it branched} along the 
real line $\an{\bf u}_{\sst\mathbb R}$. The branching occurs when $y=t=0$. 
Higher order terms deviate from the 3-plane $\an{{\bf u},{\bf v},{\bf u}
\label{bp3ex4}
\t{\bf v}}_{\sst\mathbb R}$, and make the singularity isolated.
\end{ex}

\section{Almost Calabi--Yau geometry}
\label{bp4}

{\it Calabi--Yau $m$-folds} $(M,J,\om,\Om)$ are compact complex
$m$-folds $(M,J)$ equipped with a  Ricci-flat K\"ahler metric $g$
with K\"ahler form $\om$, and a holomorphic $(m,0)$-form $\Om$ of
constant length $\ms{\Om}=2^m$. Then $\Re\Om$ is a {\it calibration}
on $(M,g)$, and the corresponding calibrated submanifolds are called
{\it special Lagrangian $m$-folds}. They are a natural generalization
of the idea of special Lagrangian submanifolds in~$\C^m$.

However, we will actually define and study special Lagrangian
submanifolds in the much larger class of {\it almost Calabi--Yau
manifolds} $(M,J,\om,\Om)$, in which $g$ is not required to be
Ricci-flat, and $\Om$ not required to have constant length. Apart
from greater generality, this has the advantage that by restricting
to a {\it generic} almost Calabi--Yau manifold one can (the author
believes) much simplify the singular behavior of the special
Lagrangian submanifolds within it.

The idea of extending special Lagrangian geometry to almost Calabi--Yau 
manifolds appears in the work of Goldstein \cite[\S 3.1]{Gold1}, Bryant 
\cite[\S 1]{Brya}, who uses the term ``special K\"ahler'' instead of 
``almost Calabi--Yau'', and the author~\cite{Joyc12}. 

\subsection{Calabi--Yau and almost Calabi--Yau manifolds}
\label{bp41}

Here is our definition of Calabi--Yau and almost Calabi--Yau manifolds.

\begin{dfn} Let $m\ge 2$. An {\it almost Calabi--Yau $m$-fold}, or
{\it ACY\/ $m$-fold}\/ for short, is a quadruple $(M,J,\om,\Om)$ 
such that $(M,J)$ is a compact $m$-dimensional complex manifold,
$\om$ is the K\"ahler form of a K\"ahler metric $g$ on $M$, and
$\Om$ is a non-vanishing holomorphic $(m,0)$-form on~$M$.

We call $(M,J,\om,\Om)$ a {\it Calabi--Yau $m$-fold}, or {\it CY\/ 
$m$-fold}\/ for short, if in addition $\om$ and $\Om$ satisfy
\e
\om^m/m!=(-1)^{m(m-1)/2}(i/2)^m\Om\w\bar\Om.
\label{bp4eq1}
\e
Furthermore, $g$ is Ricci-flat and its holonomy group is a
\label{bp4def1}
subgroup of~$\SU(m)$.
\end{dfn}

This is not the usual definition of a Calabi--Yau manifold, but
is essentially equivalent to it. Using Yau's proof of the Calabi
Conjecture \cite{Yau}, \cite[\S 5]{Joyc1} one can prove:

\begin{thm} Let\/ $(M,J)$ be a compact complex manifold with trivial
canonical bundle $K_M$, admitting K\"ahler metrics. Then in each
K\"ahler class on $M$ there is a unique Ricci-flat K\"ahler metric
$g$, with K\"ahler form $\om$. Given such $g$ and\/ $\om$, there
exists a holomorphic section $\Om$ of\/ $K_M$, unique up to change
of phase $\Om\mapsto{\rm e}^{i\th}\Om$, such that\/ $(M,J,\om,\Om)$
\label{bp4thm1}
is a Calabi--Yau manifold.
\end{thm}

Thus, to find examples of Calabi--Yau manifolds all one needs is
complex manifolds $(M,J)$ satisfying certain essentially topological
conditions. Using algebraic geometry one can construct very large
numbers of such complex manifolds, particularly in complex dimension
3, and thus Calabi--Yau manifolds are very abundant. For a review of
such constructions, and other general properties of Calabi--Yau
manifolds, see~\cite[\S 6]{Joyc1}.

\subsection{SL $m$-folds in almost Calabi--Yau $m$-folds}
\label{bp42}

Next, we define {\it special Lagrangian $m$-folds} in almost
Calabi--Yau $m$-folds.

\begin{dfn} Let $(M,J,\om,\Om)$ be an almost Calabi--Yau $m$-fold,
and $N$ a real $m$-dimensional submanifold of $M$. We call $N$ a
{\it special Lagrangian submanifold}, or {\it SL $m$-fold} for
short, if $\om\vert_N\equiv\Im\Om\vert_N\equiv 0$. It easily
follows that $\Re\Om\vert_N$ is a nonvanishing $m$-form on $N$.
Thus $N$ is orientable, with a unique orientation in which
\label{bp4def2}
$\Re\Om\vert_N$ is positive.
\end{dfn}

Let $(M,J,\om,\Om)$ be a Calabi--Yau $m$-fold, with metric $g$.
Then equation \eq{bp4eq1} ensures that for each $x\in M$ there
exists an isomorphism $T_xM\cong\C^m$ that identifies $g_x,\om_x$
and $\Om_x$ with the flat versions $g',\om',\Om'$ on $\C^m$ in
\eq{bp2eq1}. From Proposition \ref{bp2prop2} we then deduce:

\begin{prop} Let\/ $(M,J,\om,\Om)$ be a Calabi--Yau $m$-fold,
with metric $g$, and\/ $N$ a real\/ $m$-submanifold of\/ $M$.
Then $N$ is special Lagrangian, with the natural orientation,
\label{bp4prop1}
if and only if it is calibrated with respect to~$\Re\Om$.
\end{prop}

Thus, in the Calabi--Yau case Definition \ref{bp4def2} is
equivalent to the conventional definition of special Lagrangian
$m$-folds in Calabi--Yau $m$-folds, which is that they should be
calibrated with respect to $\Re\Om$, as in Definition \ref{bp2def2}.
In the almost Calabi--Yau case, we can still interpret SL $m$-folds
as calibrated submanifolds, but with respect to a conformally rescaled
metric $\ti g$. We explain how in the next proposition, which is
easily proved using Proposition~\ref{bp2prop2}.

\begin{prop} Let\/ $(M,J,\om,\Om)$ be an almost Calabi--Yau $m$-fold
with metric $g$, define $f:M\ra(0,\iy)$ by $f^{2m}\om^m/m!=(-1)^{m(m-1)/2}
(i/2)^m\Om\w\bar\Om$, and let\/ $\ti g$ be the conformally equivalent 
metric $f^2g$ on $M$. Then $\Re\Om$ is a calibration on the Riemannian 
manifold\/~$(M,\ti g)$.

A real\/ $m$-submanifold\/ $N$ in $M$ is special Lagrangian in
$(M,J,\om,\Om)$ if and only if it admits an orientation for which 
it is calibrated with respect to $\Re\Om$ in $(M,\ti g)$. In 
particular, special Lagrangian $m$-folds in $M$ are minimal 
\label{bp4prop2}
in~$(M,\ti g)$.
\end{prop}

Thus, we can give an equivalent definition of SL $m$-folds in terms
of calibrated geometry. Nonetheless, in the author's view the definition
of SL $m$-folds in terms of the vanishing of closed forms is more
fundamental than the definition in terms of calibrated geometry, at
least in the almost Calabi--Yau case, and so should be taken as the
primary definition.

One important reason for considering SL $m$-folds in almost Calabi--Yau
rather than Calabi--Yau $m$-folds is that they have much stronger
{\it genericness properties}. There are many situations in geometry
in which one uses a genericity assumption to control singular behavior.

For instance, pseudo-holomorphic curves in an arbitrary almost complex 
manifold may have bad singularities, but the possible singularities
in a generic almost complex manifold are much simpler. In the same 
way, it is reasonable to hope that in a {\it generic} Calabi--Yau 
$m$-fold, compact SL $m$-folds may have better singular behavior 
than in an arbitrary Calabi--Yau $m$-fold. 

But because Calabi--Yau manifolds come in only finite-dimensional 
families, choosing a generic Calabi--Yau structure is a fairly weak
assumption, and probably will not help very much. However, almost
Calabi--Yau manifolds come in {\it infinite-dimensional}\/ families,
so choosing a generic almost Calabi--Yau structure is a much more 
powerful thing to do, and will probably simplify the singular
behavior of compact SL $m$-folds considerably. We will return
to this idea in~Section~\ref{bp6}.

\section{Compact SL $m$-folds in ACY $m$-folds}
\label{bp5}

In this section we shall discuss {\it compact} special Lagrangian 
submanifolds in almost Calabi--Yau manifolds. Here are three
important questions which motivate work in this area.

\begin{itemize}
\item[\bf 1.] Let $N$ be a compact special
Lagrangian $m$-fold in a fixed almost Calabi--Yau $m$-fold
$(M,J,\om,\Om)$. Let ${\cal M}_\sN$ be the moduli space of 
{\it special Lagrangian deformations} of $N$, that is, 
the connected component of the set of special Lagrangian 
$m$-folds containing $N$. What can we say about ${\cal M}_\sN$? 
For instance, is it a smooth manifold, and of what dimension?
\smallskip

\item[\bf 2.] Let $\bigl\{(M,J_t,\om_t,\Om_t):t\in(-\ep,\ep)\bigr\}$ 
be a smooth 1-parameter family of almost Calabi--Yau $m$-folds. Suppose $N_0$ 
is an SL $m$-fold in $(M,J_0,\om_0,\Om_0)$. Under what conditions can we 
extend $N_0$ to a smooth family of special Lagrangian $m$-folds 
$N_t$ in $(M,J_t,\om_t,\Om_t)$ for~$t\in(-\ep,\ep)$?
\smallskip

\item[\bf 3.] In general the moduli space ${\cal M}_\sN$ in
Question 1 will be noncompact. Can we enlarge ${\cal M}_\sN$ to
a compact space $\,\,\ov{\!\!\cal M}_\sN$ by adding a ``boundary''
consisting of {\it singular} special Lagrangian $m$-folds?
If so, what is the nature of the singularities that develop?
\end{itemize}

Briefly, these questions concern the {\it deformations} of
special Lagrangian $m$-folds, {\it obstructions} to their
existence, and their {\it singularities} respectively. The
local answers to Questions 1 and 2 are well understood, and 
we shall discuss them in this section. Question 3 is the
subject of~Sections~\ref{bp6}--\ref{bp7}.

\subsection{Deformations of compact special Lagrangian $m$-folds}
\label{bp51}

The deformation theory of compact SL $m$-folds $N$ was studied by
McLean \cite{McLe}, who proved the following result in the Calabi--Yau
case. Because McLean's proof only relies on the fact that $\om\vert_N
\equiv\Im\Om\vert_N\equiv 0$, it also applies equally well to SL
$m$-folds in almost Calabi--Yau $m$-folds.

\begin{thm} Let\/ $(M,J,\om,\Om)$ be an almost Calabi--Yau $m$-fold, and\/ 
$N$ a compact special Lagrangian $m$-fold in $M$. Then the moduli 
space ${\cal M}_\sN$ of special Lagrangian deformations of\/ $N$ is 
\label{bp5thm1}
a smooth manifold of dimension $b^1(N)$, the first Betti number of\/~$N$.
\end{thm}

\begin{proof}[Sketch of proof]
Suppose for simplicity that $N$ is an
embedded submanifold. There is a natural orthogonal decomposition 
$TM\vert_N=TN\op\nu$, where $\nu\ra N$ is the {\it normal bundle} 
of $N$ in $M$. As $N$ is Lagrangian, the complex structure 
$J:TM\ra TM$ gives an isomorphism $J:\nu\ra TN$. But the metric $g$ 
gives an isomorphism $TN\cong T^*N$. Composing these two gives an
isomorphism~$\nu\cong T^*N$.

Let $T$ be a small {\it tubular neighborhood} of $N$ in $M$. Then 
we can identify $T$ with a neighborhood of the zero section in $\nu$.
Using the isomorphism $\nu\cong T^*N$, we have an identification
between $T$ and a neighborhood of the zero section in $T^*N$. This
can be chosen to identify the K\"ahler form $\om$ on $T$ with the natural
symplectic structure on $T^*N$. Let $\pi:T\ra N$ be the obvious projection.

Under this identification, submanifolds $N'$ in $T\subset M$ which 
are $C^1$ close to $N$ are identified with the graphs of small smooth 
sections $\al$ of $T^*N$. That is, submanifolds $N'$ of $M$ close to
$N$ are identified with 1-{\it forms} $\al$ on $N$. We need to know: 
which 1-forms $\al$ are identified with {\it special Lagrangian}
submanifolds~$N'$?

Well, $N'$ is special Lagrangian if $\om\vert_{N'}\equiv
\Im\Om\vert_{N'}\equiv 0$. Now $\pi\vert_{N'}:N'\ra N$ is a
diffeomorphism, so we can push $\om\vert_{N'}$ and
$\Im\Om\vert_{N'}$ down to $N$, and regard them as functions 
of $\al$. Calculation shows~that
\begin{equation*}
\pi_*\bigl(\om\vert_{N'}\bigr)=\d\al
\quad\text{and}\quad
\pi_*\bigl(\Im\Om\vert_{N'}\bigr)=F(\al,\nabla\al),
\end{equation*}
where $F$ is a nonlinear function of its arguments. Thus, the moduli 
space ${\cal M}_\sN$ is locally isomorphic to the set of small 1-forms 
$\al$ on $N$ such that $\d\al\equiv 0$ and $F(\al,\nabla\al)\equiv 0$.

Now it turns out that $F$ satisfies $F(\al,\nabla\al)\approx \d(*\al)$ 
when $\al$ is small. Therefore ${\cal M}_\sN$ is locally approximately 
isomorphic to the vector space of 1-forms $\al$ with $\d\al=\d(*\al)=0$.
But by Hodge theory, this is isomorphic to the de Rham cohomology
group $H^1(N,\R)$, and is a manifold with dimension~$b^1(N)$.

To carry out this last step rigorously requires some technical
machinery: one must work with certain {\it Banach spaces} of 
sections of $T^*N$, $\La^2T^*N$ and $\La^mT^*N$, use {\it elliptic 
regularity results} to prove that the map $\al\mapsto\bigl(\d\al,
F(\al,\nabla\al)\bigr)$ has {\it closed image} in these Banach spaces,
and then use the {\it Implicit Function Theorem for Banach spaces}
to show that the kernel of the map is what we expect.
\end{proof}

\subsection{Obstructions to the existence of compact SL $m$-folds}
\label{bp52}

Next we address Question 2 above. First, observe that if
$(M,J,\om,\Om)$ is an almost Calabi--Yau $m$-fold and $N$ a
compact SL $m$-fold in $M$ then $\om\vert_N\equiv\Im\Om\vert_N
\equiv 0$, and thus $[\om\vert_N]$ and $[\Im\Om\vert_N]$ are
zero in $H^2(N,\R)$ and $H^m(N,\R)$. But $[\om\vert_N]$ and
$[\Im\Om\vert_N]$ are unchanged under continuous variations
of $N$ in $M$. So we deduce:

\begin{lem} Let\/ $(M,J,\om,\Om)$ be an almost Calabi--Yau $m$-fold,
and\/ $N$ a compact real\/ $m$-submanifold in $M$. Then a necessary
condition for $N$ to be isotopic to a special Lagrangian submanifold\/
$N'$ in $M$ is that\/ $[\om\vert_N]=0$ in $H^2(N,\R)$ and\/
\label{bp5lem}
$[\Im\Om\vert_N]=0$ in~$H^m(N,\R)$.
\end{lem}

This gives a simple, necessary topological condition for an isotopy
class of $m$-submanifolds in an almost Calabi--Yau $m$-fold to
contain an SL $m$-fold. Our next result, following from Marshall
\cite[Th.~3.2.9]{Mars}, shows that locally, this is also a {\it
sufficient}\/ condition for an SL $m$-fold to persist under
deformations of the almost Calabi--Yau structure.

\begin{thm} Let\/ $\bigl\{(M,J_t,\om_t,\Om_t):t\in(-\ep,\ep)\bigr\}$
be a smooth family of almost Calabi--Yau $m$-folds. Let\/ $N_0$ be a
compact SL\/ $m$-fold in $(M,J_0,\om_0,\Om_0)$, and suppose that\/
$[\om_t\vert_{N_0}]=0$ in $H^2(N_0,\R)$ and\/ $[\Im\Om_t\vert_{N_0}]=0$
in $H^m(N_0,\R)$ for all\/ $t\in(-\ep,\ep)$. Then $N_0$ extends to a
smooth family $\bigl\{N_t:t\in(-\de,\de)\bigr\}$, where $0<\de\le\ep$
and\/ $N_t$ is a compact SL\/ $m$-fold in~$(M,J_t,\om_t,\Om_t)$.
\label{bp5thm2}
\end{thm}

This is proved using similar techniques to Theorem
\ref{bp5thm1}, though McLean did not prove it. Note that the
condition $[\Im\Om_t\vert_{N_0}]=0$ for all $t$ can be satisfied
by choosing the phases of the $\Om_t$ appropriately, and if the 
image of $H_2(N,\Z)$ in $H_2(M,\R)$ is zero, then the condition 
$[\om\vert_N]=0$ holds automatically. 

Thus, the obstructions $[\om_t\vert_{N_0}]=[\Im\Om_t\vert_{N_0}]=0$ 
in Theorem \ref{bp5thm2} are actually fairly mild restrictions, and 
special Lagrangian $m$-folds should be thought of as pretty stable 
under small deformations of the almost Calabi--Yau structure.
\medskip

\noindent{\bf Remark.} The deformation and obstruction theory of 
compact special Lagrangian $m$-folds are {\it extremely well-behaved}
compared to many other moduli space problems in differential geometry.
In other geometric problems (such as the deformations of complex 
structures on a complex manifold, or pseudo-holomorphic curves in 
an almost complex manifold, or instantons on a Riemannian 4-manifold, 
and so on), the deformation theory often has the following general 
structure.

There are vector bundles $E,F$ over a compact manifold $M$, and an
elliptic operator $P:C^\iy(E)\ra C^\iy(F)$, usually first-order. The 
kernel $\Ker P$ is the set of {\it infinitesimal deformations}, and 
the cokernel $\Coker P$ the set of {\it obstructions}. The actual 
moduli space $\cal M$ is locally the zeros of a nonlinear 
map~$\Psi:\Ker P\ra\Coker P$.

In a {\it generic} case, $\Coker P=0$, and then the moduli space
$\cal M$ is locally isomorphic to $\Ker P$, and so is locally a 
manifold with dimension $\ind(P)$. However, in nongeneric situations 
$\Coker P$ may be nonzero, and then the moduli space $\cal M$ may be 
nonsingular, or have an unexpected dimension.

However, SL $m$-folds do not follow this pattern. Instead, the
obstructions are {\it topologically determined}, and the moduli
space is {\it always} smooth, with dimension given by a topological
formula. This should be regarded as a minor mathematical miracle.

\subsection{Natural coordinates on the moduli space ${\cal M}_\sN$}
\label{bp53}

Let $N$ be a compact SL $m$-fold in an almost Calabi--Yau $m$-fold 
$(M,J,\om,\Om)$. Theorem \ref{bp5thm1} shows that the moduli space 
${\cal M}_\sN$ has dimension $b^1(N)$. By Poincar\'e duality 
$b^1(N)=b^{m-1}(N)$. Thus ${\cal M}_\sN$ has the same dimension as 
the de Rham cohomology groups $H^1(M,\R)$ and~$H^{m-1}(M,\R)$. 

We shall construct natural local diffeomorphisms $\Phi$ from
${\cal M}_\sN$ to $H^1(N,\R)$, and $\Psi$ from ${\cal M}_\sN$ to
$H^{m-1}(N,\R)$. These induce two natural {\it affine structures} 
on ${\cal M}_\sN$, and can be thought of as two {\it natural 
coordinate systems} on ${\cal M}_\sN$. The material of this
section can be found in Hitchin~\cite[\S 4]{Hitc}.

Here is how to define $\Phi$ and $\Psi$. Let $U$ be a connected
and simply-connected open neighborhood of $N$ in ${\cal M}_\sN$. 
We will construct smooth maps $\Phi:U\ra H^1(N,\R)$ and
$\Psi:U\ra H^{m-1}(N,\R)$ with $\Phi(N)=\Psi(N)=0$, which are 
local diffeomorphisms.

Let $N'\in U$. Then as $U$ is connected, there exists a smooth
path $\ga:[0,1]\ra U$ with $\ga(0)=N$ and $\ga(1)=N'$, and as
$U$ is simply-connected, $\ga$ is unique up to isotopy. Now 
$\ga$ parametrizes a family of submanifolds of $M$ diffeomorphic
to $N$, which we can lift to a smooth map $\Ga:N\t[0,1]\ra M$
with~$\Ga(N\t\{t\})=\ga(t)$.

Consider the 2-form $\Ga^*(\om)$ on $N\t[0,1]$. As each fiber
$\ga(t)$ is Lagrangian, we have $\Ga^*(\om)\vert_{N\t\{t\}}\equiv 0$
for each $t\in[0,1]$. Therefore we may write $\Ga^*(\om)=\al_t\w\d t$,
where $\al_t$ is a closed 1-form on $N$ for $t\in[0,1]$. Define
\begin{equation*}
\Phi(N')=\bigl[\ts\int_0^1\al_t\,\d t\bigr]\in H^1(N,\R).
\end{equation*}
That is, we
integrate the 1-forms $\al_t$ with respect to $t$ to get a closed 
1-form $\int_0^1\al_t\,\d t$, and then take its cohomology class. 

Similarly, write $\Ga^*(\Im\Om)=\be_t\w\d t$, where $\be_t$ is a 
closed $(m\!-\!1)$-form on $N$ for $t\in[0,1]$, and define
$\Psi(N')=\bigl[\int_0^1\be_t\,\d t\bigr]\in H^{m-1}(N,\R)$.
Then $\Phi$ and $\Psi$ are independent of choices made in the 
construction (exercise). We need to restrict to a simply-connected 
subset $U$ of ${\cal M}_\sN$ so that $\ga$ is unique up to isotopy. 
Alternatively, one can define $\Phi$ and $\Psi$ on the universal 
cover $\,\,\widetilde{\!\!\cal M}_\sN$ of~${\cal M}_\sN$.

\section{Singularities of special Lagrangian $m$-folds}
\label{bp6}

Now we move on to Question 3 of Section~\ref{bp5}, and discuss the
{\it singularities} of special Lagrangian $m$-folds. We can
divide it into two sub-questions:
\begin{itemize}
\item[\bf 3(a)] What kinds of singularities are possible in
singular special Lagrangian $m$-folds, and what do they look like?
\item[\bf 3(b)] How can singular SL $m$-folds arise as limits
of nonsingular SL $m$-folds, and what does the limiting
behavior look like near the singularities?
\end{itemize}

These questions are addressed in the author's series of papers
\cite{Joyc13,Joyc14,Joyc15,Joyc16,Joyc17} for {\it isolated
conical singularities} of SL $m$-folds, that is, singularities
locally modelled on an SL cone $C$ in $\C^m$ with an isolated
singularity at 0. We now explain the principal results. Readers
of the series are advised to begin with the final paper
\cite{Joyc17}, which surveys the others.

\subsection{Special Lagrangian cones}
\label{bp61}

We define {\it SL cones}, and some notation.

\begin{dfn} A (singular) SL $m$-fold $C$ in $\C^m$ is called a
{\it cone} if $C=tC$ for all $t>0$, where $tC=\{t\,{\bf x}:{\bf x}
\in C\}$. Let $C$ be a closed SL cone in $\C^m$ with an isolated
singularity at 0. Then $\Si=C\cap{\cal S}^{2m-1}$ is a compact,
nonsingular $(m\!-\!1)$-submanifold of ${\cal S}^{2m-1}$, not
necessarily connected. Let $g_\sSi$ be the restriction
of $g'$ to $\Si$, where $g'$ is as in~\eq{bp2eq1}.

Set $C'=C\sm\{0\}$. Define $\iota:\Si\t(0,\iy)\ra\C^m$ by
$\iota(\si,r)=r\si$. Then $\iota$ has image $C'$. By an abuse
of notation, {\it identify} $C'$ with $\Si\t(0,\iy)$ using
$\iota$. The {\it cone metric} on $C'\cong\Si\t(0,\iy)$ is
$g'=\iota^*(g')=\d r^2+r^2g_\sSi$. For $\al\in\R$, we say
that a function $u:C'\ra\R$ is {\it homogeneous of order}
$\al$ if $u\circ t\equiv t^\al u$ for all $t>0$. Equivalently,
\label{bp6def1}
$u$ is homogeneous of order $\al$ if $u(\si,r)\equiv r^\al
v(\si)$ for some function~$v:\Si\ra\R$.
\end{dfn}

In \cite[Lem.~2.3]{Joyc13} we study {\it homogeneous harmonic
functions} on~$C'$.

\begin{lem} In the situation of Definition \ref{bp6def1},
let\/ $u(\si,r)\equiv r^\al v(\si)$ be a homogeneous function
of order $\al$ on $C'=\Si\t(0,\iy)$, for $v\in C^2(\Si)$. Then
\begin{equation*}
\De u(\si,r)=r^{\al-2}\bigl(\De_\sSi v-\al(\al+m-2)v\bigr),
\end{equation*}
where $\De$, $\De_\sSi$ are the Laplacians on $(C',g')$
and\/ $(\Si,g_\sSi)$. Hence, $u$ is harmonic on $C'$
if and only if\/ $v$ is an eigenfunction of\/ $\De_\sSi$
\label{bp6lem}
with eigenvalue~$\al(\al+m-2)$.
\end{lem}

Following \cite[Def.~2.5]{Joyc13}, we define:

\begin{dfn} In the situation of Definition \ref{bp6def1},
suppose $m>2$ and define
\e
\D_\sSi=\bigl\{\al\in\R:\text{$\al(\al+m-2)$ is
an eigenvalue of $\De_\sSi$}\bigr\}.
\label{bp6eq1}
\e
Then $\D_\sSi$ is a countable, discrete subset of
$\R$. By Lemma \ref{bp6lem}, an equivalent definition is that
$\D_\sSi$ is the set of $\al\in\R$ for which there
exists a nonzero homogeneous harmonic function $u$ of order
$\al$ on~$C'$.

Define $m_\sSi:\D_\sSi\ra\N$ by taking
$m_\sSi(\al)$ to be the multiplicity of the eigenvalue
$\al(\al+m-2)$ of $\De_\sSi$, or equivalently the
dimension of the vector space of homogeneous harmonic
functions $u$ of order $\al$ on $C'$. Define
$N_\sSi:\R\ra\Z$ by
\e
N_\sSi(\de)=
-\sum_{\!\!\!\!\al\in\D_\sSi\cap(\de,0)\!\!\!\!}m_\sSi(\al)
\;\>\text{if $\de<0$, and}\;\>
N_\sSi(\de)=
\sum_{\!\!\!\!\al\in\D_\sSi\cap[0,\de]\!\!\!\!}m_\sSi(\al)
\;\>\text{if $\de\ge 0$.}
\label{bp6eq2}
\e
Then $N_\sSi$ is monotone increasing and upper semicontinuous,
and is discontinuous exactly on $\D_\sSi$, increasing by
$m_\sSi(\al)$ at each $\al\in\D_\sSi$. As the
eigenvalues of $\De_\sSi$ are nonnegative, we see that
\label{bp6def2}
$\D_\sSi\cap(2-m,0)=\emptyset$ and $N_\sSi\equiv 0$
on~$(2-m,0)$.
\end{dfn}

We define the {\it stability index} of $C$, and {\it stable}
cones~\cite[Def.~3.6]{Joyc14}.

\begin{dfn} Let $C$ be an SL cone in $\C^m$ for $m>2$ with an
isolated singularity at 0, let $G$ be the Lie subgroup of\/
$\SU(m)$ preserving $C$, and use the notation of Definitions
\ref{bp6def1} and \ref{bp6def2}. Then \cite[eq.~(8)]{Joyc14}
shows that
\e
m_\sSi(0)=b^0(\Si),\quad
m_\sSi(1)\ge 2m \quad\text{and}\quad
m_\sSi(2)\ge m^2-1-\dim G.
\label{bp6eq3}
\e

Define the {\it stability index} $\sind(C)$ to be
\e
\sind(C)=N_\sSi(2)-b^0(\Si)-m^2-2m+1+\dim G.
\label{bp6eq4}
\e
Then $\sind(C)\ge 0$ by \eq{bp6eq3}, as $N_\sSi(2)\ge
m_\sSi(0)+m_\sSi(1)+m_\sSi(2)$ by \eq{bp6eq2}.
\label{bp6def3}
We call $C$ {\it stable} if~$\sind(C)=0$.
\end{dfn}

Here is the point of these definitions. By the Principle in
Section \ref{bp23}, homogeneous harmonic functions $v$ of
order $\al$ on $C'$ correspond to infinitesimal deformations
$\d v$ of $C'$ as an SL $m$-fold in $\C^m$, which grow like
$O(r^{\al-1})$. Hence, $N_\sSi(\la)$ is effectively the dimension
of a space of {\it infinitesimal deformations} of $C'$ as an SL
$m$-fold, which grow like $O(r^{\al-1})$ for $\al\in[0,\la]$,
when~$\la\ge 0$.

For $\la=2$ this space of harmonic functions, or infinitesimal
deformations of $C'$, contains some from obvious geometrical
sources: locally constant functions on $C'$, and infinitesimal
deformations of $C'$ from translations in $\C^m$, or $\su(m)$
rotations. We get $\sind(C)$ by subtracting off these obvious
geometrical deformations from $N_\sSi(2)$. Hence, $\sind(C)$ is
the dimension of a space of {\it excess infinitesimal deformations}
of $C'$ as an SL $m$-fold, with growth between $O(r^{-1})$ and
$O(r)$, which do not arise from infinitesimal automorphisms
of~$\C^m$.

We shall see in Section~\ref{bp62} that $\sind(C)$ is the dimension
of an {\it obstruction space} to deforming an SL $m$-fold $X$
with a conical singularity with cone $C$, and that if $C$ is
{\it stable} then the deformation theory of $X$ simplifies.

\subsection{Special Lagrangian $m$-folds with conical singularities}
\label{bp62}

Now we can define {\it conical singularities} of SL $m$-folds,
following~\cite[Def.~3.6]{Joyc13}.

\begin{dfn} Let $(M,J,\om,\Om)$ be an almost Calabi--Yau $m$-fold
for $m>2$. Suppose $X$ is a compact singular SL $m$-fold in $M$
with singularities at distinct points $x_1,\ldots,x_n\in X$, and
no other singularities.

Fix isomorphisms $\up_i:\C^m\ra T_{x_i}M$ for $i=1,\ldots,n$
such that $\up_i^*(\om)=\om'$ and $\up_i^*(\Om)=a_i\Om'$, where
$\om',\Om'$ are as in \eq{bp2eq1} and $a_1,\ldots,a_n>0$. Let
$C_1,\ldots,C_n$ be SL cones in $\C^m$ with isolated singularities
at 0. For $i=1,\ldots,n$ let $\Si_i=C_i\cap{\cal S}^{2m-1}$, and
let $\mu_i\in(2,3)$ with
\e
(2,\mu_i]\cap\D_\sSii=\emptyset,
\quad\text{where $\D_\sSii$ is defined in \eq{bp6eq1}.}
\label{bp6eq5}
\e
Then we say that
$X$ has a {\it conical singularity} or {\it conical singular
point} at $x_i$, with {\it rate} $\mu_i$ and {\it cone} $C_i$
for $i=1,\ldots,n$, if the following holds.

By Darboux Theorem there exist embeddings $\Up_i:B_R\ra M$ for
$i=1,\ldots,n$ satisfying $\Up_i(0)=x_i$, $\d\Up_i\vert_0=\up_i$
and $\Up_i^*(\om)=\om'$, where $B_R$ is the open ball of radius $R$
about 0 in $\C^m$ for some small $R>0$. Define $\iota_i:\Si_i\t(0,R)
\ra B_R$ by $\iota_i(\si,r)=r\si$ for~$i=1,\ldots,n$.

Define $X'=X\sm\{x_1,\ldots,x_n\}$. Then there should exist a
compact subset $K\subset X'$ such that $X'\sm K$ is a union of
open sets $S_1,\ldots,S_n$ with $S_i\subset\Up_i(B_R)$, whose
closures $\bar S_1,\ldots,\bar S_n$ are disjoint in $X$. For
$i=1,\ldots,n$ and some $R'\in(0,R]$ there should exist a smooth
$\phi_i:\Si_i\t(0,R')\ra B_R$ such that $\Up_i\circ\phi_i:\Si_i
\t(0,R')\ra M$ is a diffeomorphism $\Si_i\t(0,R')\ra S_i$, and
\e
\bmd{\nabla^k(\phi_i-\iota_i)}=O(r^{\mu_i-1-k})
\quad\text{as $r\ra 0$ for $k=0,1$.}
\label{bp6eq6}
\e
Here $\nabla$ is the Levi-Civita connection of the cone metric
$\iota_i^*(g')$ on $\Si_i\t(0,R')$, $\md{\,.\,}$ is computed
using $\iota_i^*(g')$. If the cones $C_1,\ldots,C_n$ are
{\it stable} in the sense of Definition \ref{bp6def3}, then
\label{bp6def4}
we say that $X$ has {\it stable conical singularities}.
\end{dfn}

We show in \cite[Th.s 4.4 \& 5.5]{Joyc13} that if \eq{bp6eq6}
holds for $k=0,1$ and some $\mu_i$ satisfying \eq{bp6eq5}, then
we can choose a natural $\phi_i$ for which \eq{bp6eq6} holds
for {\it all\/} $k\ge 0$, and for {\it all\/} rates $\mu_i$
satisfying \eq{bp6eq5}. Thus the number of derivatives
required in \eq{bp6eq6} and the choice of $\mu_i$ both make
little difference. We choose $k=0,1$ in \eq{bp6eq6}, and some
$\mu_i$ in \eq{bp6eq5}, to make the definition as weak
as possible.

Suppose we did not require \eq{bp6eq5}, and that $\al\in(2,\mu_i)
\cap\D_\sSii$. Then there exists a homogeneous harmonic
function $v$ on $C_i'$ of order $\al$. By the Principle in
Section \ref{bp23}, $\d v$ yields an infinitesimal
deformation of $C_i'$ as an SL $m$-fold in $\C^m$, growing
like $O(r^{\al-1})$. Locally this gives a way to deform $X$
into an SL $m$-fold which would {\it not\/} satisfy \eq{bp6eq6},
as $\al<\mu_i$. Effectively, $v$ acts as an {\it obstruction}
to deforming $X$ through SL $m$-folds satisfying Definition
\ref{bp6def4}. So the point of \eq{bp6eq5} is to reduce to
a minimum the obstructions to existence and deformation of
SL $m$-folds with isolated conical singularities.

In \cite{Joyc14} we study the {\it deformation theory} of
compact SL $m$-folds with conical singularities, generalizing
Theorem \ref{bp5thm1} in the nonsingular case. Following
\cite[Def.~5.4]{Joyc14}, we define the space $\M_\sX$ of
compact SL $m$-folds $\hat X$ in $M$ with conical singularities
deforming a fixed SL $m$-fold $X$ with conical singularities.

\begin{dfn} Let $(M,J,\om,\Om)$ be an almost Calabi--Yau
$m$-fold and $X$ a compact SL $m$-fold in $M$ with conical
singularities at $x_1,\ldots,x_n$ with identifications
$\up_i:\C^m\ra T_{x_i}M$ and cones $C_1,\ldots,C_n$. Define
the {\it moduli space} $\M_\sX$ {\it of deformations of\/}
$X$ to be the set of $\hat X$ such that
\label{bp6def5}
\begin{itemize}
\setlength{\parsep}{0pt}
\setlength{\itemsep}{0pt}
\item[(i)] $\hat X$ is a compact SL $m$-fold in $M$ with
conical singularities at $\hat x_1,\ldots,\hat x_n$ with
cones $C_1,\ldots,C_n$, for some $\hat x_i$ and
identifications~$\hat\up_i:\C^m\ra T_{\smash{\hat x_i}}M$.
\item[(ii)] There exists a homeomorphism $\hat\iota:X\ra\hat X$
with $\hat\iota(x_i)=\hat x_i$ for $i=1,\ldots,n$ such that
$\hat\iota\vert_{X'}:X'\ra\hat X'$ is a diffeomorphism and
$\hat\iota$ and $\iota$ are isotopic as continuous maps
$X\ra M$, where $\iota:X\ra M$ is the inclusion.
\end{itemize}
\end{dfn}

In \cite[Def.~5.6]{Joyc14} we define a {\it topology} on
$\M_\sX$, and explain why it is the natural choice. We
will not repeat the complicated definition here. In
\cite[Th.~6.10]{Joyc14} we describe $\M_\sX$ near $X$, in terms
of a smooth map $\Phi$ between the {\it infinitesimal deformation
space} $\I_\sXp$ and the {\it obstruction space}~$\O_\sXp$.

\begin{thm} Suppose $(M,J,\om,\Om)$ is an almost Calabi--Yau
$m$-fold and\/ $X$ a compact SL\/ $m$-fold in $M$ with conical
singularities at\/ $x_1,\ldots,x_n$ and cones $C_1,\ldots,C_n$.
Let\/ $\M_\sX$ be the moduli space of deformations
of\/ $X$ as an SL\/ $m$-fold with conical singularities in $M$,
as in Definition \ref{bp6def5}. Set\/~$X'=X\sm\{x_1,\ldots,x_n\}$.

Then there exist natural finite-dimensional vector spaces
$\I_\sXp$, $\O_\sXp$ such that\/ $\I_\sXp$ is the image of\/
$H^1_{\rm cs}(X',\R)$ in $H^1(X',\R)$ and\/ $\dim\O_\sXp=
\sum_{i=1}^n\sind(C_i)$, where $\sind(C_i)$ is the stability
index of Definition \ref{bp6def3}. There exists an open
neighbourhood\/ $U$ of\/ $0$ in $\I_\sXp$, a smooth
map $\Phi:U\ra\O_\sXp$ with\/ $\Phi(0)=0$, and a map
$\Xi:\{u\in U:\Phi(u)=0\}\ra\M_\sX$ with\/ $\Xi(0)=X$ which is
\label{bp6thm1}
a homeomorphism with an open neighbourhood of\/ $X$ in~$\M_\sX$.
\end{thm}

If the $C_i$ are {\it stable} then $\O_\sXp=\{0\}$ and we
deduce~\cite[Cor.~6.11]{Joyc14}:

\begin{cor} Suppose $(M,J,\om,\Om)$ is an almost Calabi--Yau
$m$-fold and\/ $X$ a compact SL\/ $m$-fold in $M$ with stable
conical singularities, and let\/ $\M_\sX$ and\/ $\I_\sXp$ be
as in Theorem \ref{bp6thm1}. Then $\M_\sX$ is a smooth manifold
\label{bp6cor1}
of dimension~$\dim\I_\sXp$.
\end{cor}

This has clear similarities with Theorem \ref{bp5thm1}. Here
is another simple condition for $\M_\sX$ to be a manifold
near $X$, \cite[Def.~6.12]{Joyc14}.

\begin{dfn} Let $(M,J,\om,\Om)$ be an almost Calabi--Yau
$m$-fold and $X$ a compact SL $m$-fold in $M$ with conical
singularities, and let $\I_\sXp,\O_\sXp,U$ and $\Phi$ be as
in Theorem \ref{bp6thm1}. We call $X$ {\it transverse} if
\label{bp6def6}
the linear map $\d\Phi\vert_0:\I_\sXp\ra\O_\sXp$ is surjective.
\end{dfn}

If $X$ is transverse then $\{u\in U:\Phi(u)=0\}$ is a manifold
near 0, so Theorem \ref{bp6thm1} yields~\cite[Cor.~6.13]{Joyc14}:

\begin{cor} Suppose $(M,J,\om,\Om)$ is an almost Calabi--Yau $m$-fold
and\/ $X$ a transverse compact SL\/ $m$-fold in $M$ with conical
singularities, and let\/ $\M_\sX,\I_\sXp$ and\/ $\O_\sXp$ be as
in Theorem \ref{bp6thm1}. Then $\M_\sX$ is near $X$ a smooth manifold
\label{bp6cor2}
of dimension~$\dim\I_\sXp-\dim\O_\sXp$.
\end{cor}

We would like to conclude that by choosing a sufficiently
generic perturbation $\om^s$ we can make $\M_\sX^s$
smooth everywhere. This is the idea of the following
conjecture,~\cite[Conj.~9.5]{Joyc14}:

\begin{conj} Let\/ $(M,J,\om,\Om)$ be an almost Calabi--Yau
$m$-fold, $X$ a compact SL\/ $m$-fold in $M$ with conical
singularities, and\/ $\I_\sXp,\O_\sXp$ be as in Theorem
\ref{bp6thm1}. Then for a second category subset of K\"ahler
forms $\hat\om$ in the K\"ahler class of\/ $\om$, the moduli
space $\hat\M_\sX$ of compact SL\/ $m$-folds $\hat X$ with
conical singularities in $(M,J,\hat\om,\Om)$ isotopic to $X$
consists of transverse $\hat X$, and so is a smooth manifold
of dimension~$\dim\I_\sXp-\dim\O_\sXp$.
\label{bp6conj}
\end{conj}

For a partial proof of this, see \cite[Th.s~9.1 \& 9.3]{Joyc14}.
Basically, we can prove the conjecture for $\hat\om$ close
to $\om$ and $\hat X\in\hat\M_\sX$ close to $X$, or more
generally, close to a fixed compact subset of the moduli
space $\M_\sX$ in~$(M,J,\om,\Om)$.

\subsection{Asymptotically Conical SL $m$-folds}
\label{bp63}

The local models for how to {\it desingularize} compact
SL $m$-folds with isolated conical singularities are
{\it Asymptotically Conical\/} SL $m$-folds $L$ in
$\C^m$, so we discuss these briefly. Here is the
definition,~\cite[Def.~7.1]{Joyc13}.

\begin{dfn} Let $C$ be a closed SL cone in $\C^m$ with isolated
singularity at 0 for $m>2$, and let $\Si=C\cap{\cal S}^{2m-1}$,
so that $\Si$ is a compact, nonsingular $(m-1)$-manifold, not
necessarily connected. Let $g_\sSi$ be the metric on $\Si$
induced by the metric $g'$ on $\C^m$ in \eq{bp2eq1}, and $r$ the
radius function on $\C^m$. Define $\iota:\Si\t(0,\iy)\ra\C^m$ by
$\iota(\si,r)=r\si$. Then the image of $\iota$ is $C\sm\{0\}$,
and $\iota^*(g')=r^2g_\sSi+\d r^2$ is the cone metric
on~$C\sm\{0\}$.

Let $L$ be a closed, nonsingular SL $m$-fold in $\C^m$. We
call $L$ {\it Asymptotically Conical (AC)} with {\it rate}
$\la<2$ and {\it cone} $C$ if there exists a compact subset
$K\subset L$ and a diffeomorphism $\vp:\Si\t(T,\iy)\ra L\sm K$
for some $T>0$, such that
\begin{equation*}
\bmd{\nabla^k(\vp-\iota)}=O(r^{\la-1-k})
\quad\text{as $r\ra\iy$ for $k=0,1$.}
\end{equation*}
Here $\nabla,\md{\,.\,}$ are computed
\label{bp6def7}
using the cone metric~$\iota^*(g')$.
\end{dfn}

The deformation theory of Asymptotically Conical SL $m$-folds
in $\C^m$ has been studied independently by Pacini \cite{Paci}
and Marshall \cite{Mars}. Pacini's results are earlier, but
Marshall's are more complete.

\begin{dfn} Suppose $L$ is an Asymptotically Conical SL
$m$-fold in $\C^m$ with cone $C$ and rate $\la<2$, as in
Definition \ref{bp6def7}. Define the {\it moduli space
$\M_\sL^\la$ of deformations of\/ $L$ with rate} $\la$
to be the set of AC SL $m$-folds $\hat L$ in $\C^m$ with
cone $C$ and rate $\la$, such that $\hat L$ is diffeomorphic
to $L$ and isotopic to $L$ as an Asymptotically Conical
\label{bp6def8}
submanifold of $\C^m$. One can define a natural {\it topology}
on~$\M_\sL^\la$.
\end{dfn}

The following result can be deduced from Marshall
\cite[Th.~6.2.15]{Mars} and \cite[Table~5.1]{Mars}.
(See also Pacini \cite[Th.~2 \& Th.~3]{Paci}.)

\begin{thm} Let\/ $L$ be an Asymptotically Conical SL\/
$m$-fold in $\C^m$ with cone $C$ and rate $\la<2$, and
let\/ $\M_\sL^\la$ be as in Definition \ref{bp6def8}.
Set\/ $\Si=C\cap{\cal S}^{2m-1}$, and let\/ $\D_\sSi,N_\sSi$
be as in Section \ref{bp61} and\/ $b^k(L),b^k_{\rm cs}(L)$
be the Betti numbers in ordinary and compactly-supported
de Rham cohomology $H^k(L,\R)$, $H^k_{\rm cs}(L,\R)$. Then
\label{bp6thm2}
\begin{itemize}
\item[{\rm(a)}] If\/ $\la\in(0,2)\sm\D_\sSi$ then
$\M_\sL^\la$ is a manifold with
\e
\dim\M_\sL^\la=b^1(L)-b^0(L)+N_\sSi(\la).
\label{bp6eq7}
\e
Note that if\/ $0<\la<\min\bigl(\D_\sSi\cap
(0,\iy)\bigr)$ then~$N_\sSi(\la)=b^0(\Si)$.
\item[{\rm(b)}] If\/ $\la\in(2-m,0)$ then $\M_\sL^\la$
is a manifold of dimension~$b^1_{\rm cs}(L)=b^{m-1}(L)$.
\end{itemize}
\end{thm}

This is the analogue of Theorems \ref{bp5thm1} and \ref{bp6thm1}
for AC SL $m$-folds. If $\la\in(2-m,2)\sm\D_\sSi$ then the
deformation theory for $L$ with rate $\la$ is {\it unobstructed\/}
and $\M_\sL^\la$ is a {\it smooth manifold\/} with a given dimension.

\subsection{Desingularizing singular SL $m$-folds}
\label{bp64}

Suppose $(M,J,\om,\Om)$ is an almost Calabi--Yau $m$-fold, and
$X$ a compact SL $m$-fold in $M$ with conical singularities
at $x_1,\ldots,x_n$ and cones $C_1,\ldots,C_n$. In
\cite{Joyc15,Joyc16} we study {\it desingularizations} of $X$,
realizing $X$ as a limit of a family of compact, nonsingular
SL $m$-folds $\smash{\ti N^t}$ in $M$ for small~$t>0$.

Here is the basic method. Let $L_1,\ldots,L_n$ be {\it Asymptotically
Conical\/} SL $m$-folds in $\C^m$, as in Section \ref{bp63}, with $L_i$
asymptotic to the cone $C_i$ at infinity. We shrink $L_i$ by a small
factor $t>0$, and glue $tL_i$ into $X$ at $x_i$ for $i=1\ldots,n$ to
get a 1-parameter family of compact, nonsingular {\it Lagrangian}
$m$-folds $N^t$ in $(M,\om)$ for small~$t>0$.

Then we show using analysis that when $t$ is sufficiently small
we can deform $N^t$ to a compact, nonsingular {\it special\/}
Lagrangian $m$-fold $\smash{\ti N^t}$ via a small Hamiltonian
deformation. This $\smash{\ti N^t}$ depends smoothly on $t$,
and as $t\ra 0$ it converges to the singular SL $m$-fold $X$,
in the sense of currents.

Our simplest desingularization result is~\cite[Th.~6.13]{Joyc15}.

\begin{thm} Let\/ $(M,J,\om,\Om)$ be an almost Calabi--Yau
$m$-fold and\/ $X$ a compact SL\/ $m$-fold in $M$ with conical
singularities at\/ $x_1,\ldots,x_n$ and cones $C_1,\ldots,C_n$.
Let\/ $L_1,\ldots,L_n$ be Asymptotically Conical SL\/ $m$-folds
in $\C^m$ with cones $C_1,\ldots,C_n$ and rates $\la_1,\ldots,
\la_n$. Suppose $\la_i<0$ for $i=1,\ldots,n$, and\/ $X'=X\sm
\{x_1,\ldots,x_n\}$ is connected.

Then there exists $\ep>0$ and a smooth family $\bigl\{
\smash{\ti N^t}:t\in(0,\ep]\bigr\}$ of compact, nonsingular
SL\/ $m$-folds in $(M,J,\om,\Om)$, such that\/ $\smash{\ti N^t}$
is constructed by gluing $tL_i$ into $X$ at\/ $x_i$ for
$i=1,\ldots,n$. In the sense of currents, $\smash{\ti N^t}\ra X$
as~$t\ra 0$.
\label{bp6thm3}
\end{thm}

The theorem contains two {\it simplifying assumptions}: that
$\la_i<0$ for all $i$, and that $X'$ is connected. These avoid
two kinds of {\it obstructions} to desingularizing $X$ using
the $L_i$. For the first, the $L_i$ have {\it cohomological
invariants} $Y(L_i)$ in $H^1(\Si_i,\R)$ derived from the
relative cohomology class of $\om'$. If $\la_i<0$ then
$Y(L_i)=0$. But if $\la_i\ge 0$ and $Y(L_i)\ne 0$ then
there are obstructions to the existence of $N^t$ as a
{\it Lagrangian} $m$-fold. That is, we can only define
$N^t$ if the $Y(L_i)$ satisfy an equation.

For the second, if $X'$ is not connected then there is an
analytic obstruction to deforming $N^t$ to $\smash{\ti N^t}$,
because the Laplacian $\De$ on functions on $N^t$ has {\it
small eigenvalues}. Again, the $L_i$ have cohomological
invariants $Z(L_i)$ in $H^{m-1}(\Si_i,\R)$ derived
from the relative cohomology class of $\Im\Om'$, and
we can only deform $N^t$ to $\smash{\ti N^t}$ if the
$Z(L_i)$ satisfy an equation.

In the obstructed cases we prove generalizations of
Theorem \ref{bp6thm3} showing that SL desingularizations
$\smash{\ti N^t}$ exist when $Y(L_i),Z(L_i)$ satisfy
equations, and also generalize the results to {\it
families} of almost Calabi--Yau $m$-folds. As the
details are complicated we will not give them, but
we refer the reader to \cite{Joyc15,Joyc16}
and~\cite[\S 7]{Joyc17}.

\subsection{The index of singularities of SL\/ $m$-folds}
\label{bp65}

We now consider the {\it boundary\/} $\pd\M_\sN$ of a
moduli space $\M_\sN$ of SL $m$-folds.

\begin{dfn} Let $(M,J,\om,\Om)$ be an almost Calabi--Yau
$m$-fold, $N$ a compact, nonsingular SL $m$-fold in $M$,
and $\M_\sN$ the moduli space of deformations of $N$ in $M$.
Then $\M_\sN$ is a smooth manifold of dimension $b^1(N)$,
in general noncompact. We can construct a natural {\it
compactification} $\oM_\sN$ as follows.

Regard $\M_\sN$ as a moduli space of special Lagrangian
{\it integral currents} in the sense of Geometric Measure
Theory, as discussed in \cite[\S 6]{Joyc13}. Let $\oM_\sN$
be the closure of $\M_\sN$ in the space of integral currents.
As elements of $\M_\sN$ have uniformly bounded volume,
$\oM_\sN$ is {\it compact}. Define the {\it boundary}
$\pd\M_\sN$ to be $\oM_\sN\sm\M_\sN$. Then elements of
\label{bp6def9}
$\pd\M_\sN$ are {\it singular SL integral currents}.
\end{dfn}

In good cases, say if $(M,J,\om,\Om)$ is suitably generic,
it seems reasonable that $\pd\M_\sN$ should be divided into
a number of {\it strata}, each of which is a moduli space of
singular SL $m$-folds with singularities of a particular type,
and is itself a manifold with singularities. In particular, some
or all of these strata could be moduli spaces $\M_\sX$ of SL
$m$-folds with isolated conical singularities, as in Section~\ref{bp62}. 

Let $\M_\sN$ be a moduli space of compact, nonsingular SL
$m$-folds $N$ in $M$, and $\M_\sX$ a moduli space of singular
SL $m$-folds in $\pd\M_\sN$ with singularities of a particular
type, and $X\in\M_\sX$. Following \cite[\S 8.3]{Joyc17}, we
(loosely) define the {\it index} of the singularities of $X$
to be $\ind(X)=\dim\M_\sN-\dim\M_\sX$, provided $\M_\sX$ is
smooth near $X$. Note that $\ind(X)$ depends on $N$ as well as~$X$.

In \cite[Th.~8.10]{Joyc17} we use the results of
\cite{Joyc14,Joyc15,Joyc16} to compute $\ind(X)$ when $X$
is {\it transverse} with conical singularities, in the sense
of Definition \ref{bp6def6}. Here is a simplified version of
the result, where we assume that $H^1_{\rm cs}(L_i,\R)\ra
H^1(L_i,\R)$ is surjective to avoid a complicated correction
term to $\ind(X)$ related to the obstructions to defining
$N^t$ as a Lagrangian $m$-fold.

\begin{thm} Let\/ $X$ be a compact, transverse SL $m$-fold in
$(M,J,\om,\Om)$ with conical singularities at\/ $x_1,\ldots,x_n$
and cones $C_1,\ldots,C_n$. Let\/ $L_1,\ldots,L_n$ be AC SL\/
$m$-folds in $\C^m$ with cones $C_1,\ldots,C_n$, such that the
natural projection $H^1_{\rm cs}(L_i,\R)\ra H^1(L_i,\R)$ is
surjective. Construct desingularizations $N$ of\/ $X$ by gluing
AC SL\/ $m$-folds $L_1,\ldots,L_n$ in at\/ $x_1,\ldots,x_n$, as
\label{bp6thm4}
in Section \ref{bp64}. Then
\e
\ind(X)=1-b^0(X')+\ts\sum_{i=1}^nb^1_{\rm cs}(L_i)+
\sum_{i=1}^n\sind(C_i).
\label{bp6eq8}
\e
\end{thm}

If the cones $C_i$ are not {\it rigid\/}, for instance if
$C_i\sm\{0\}$ is not connected, then \eq{bp6eq8} should
be corrected, as in \cite[\S 8.3]{Joyc17}. If Conjecture
\ref{bp6conj} is true then for a {\it generic} K\"ahler
form $\om$, {\it all\/} compact SL $m$-folds $X$ with
conical singularities are transverse, and so Theorem
\ref{bp6thm4} and \cite[Th.~8.10]{Joyc17} allow us to
calculate~$\ind(X)$. 

Now singularities with {\it small index} are the most commonly
occurring, and so arguably the most interesting kinds of
singularity. Also, as $\ind(X)\le\dim\M_\sN$, for various
problems it will only be necessary to know about singularities
with index up to a certain value.

For example, in \cite{Joyc4} the author proposed to define an
invariant of almost Calabi--Yau 3-folds by counting special 
Lagrangian homology 3-spheres (which occur in 0-dimensional 
moduli spaces) in a given homology class, with a certain 
topological weight. This invariant will only be interesting 
if it is essentially conserved under deformations of the 
underlying almost Calabi--Yau 3-fold. During such a deformation, 
nonsingular SL 3-folds can develop singularities and disappear, 
or new ones appear, which might change the invariant.

To prove the invariant is conserved, we need to show that it is 
unchanged along generic 1-parameter families of almost Calabi--Yau
3-folds. The only kinds of singularities of SL homology 3-spheres
that arise in such families will have index 1. Thus, to resolve
the conjectures in \cite{Joyc4}, we only have to know about index
1 singularities of SL 3-folds in almost Calabi--Yau 3-folds.

Another problem in which the index of singularities will be important
is the {\it SYZ Conjecture}, to be discussed in Section \ref{bp7}.
This has to do with dual 3-dimensional families ${\cal F},\hat{\cal F}$
of SL 3-tori in (almost) Calabi--Yau 3-folds $M,\hat M$. If
$M,\hat M$ are generic then the only kinds of singularities
that can occur at the boundaries of ${\cal F},\hat{\cal F}$ are of
index 1, 2 or 3. So, to study the SYZ Conjecture in the generic
case, we only have to know about singularities of SL 3-folds with
index 1, 2 and~3.

\section{The SYZ Conjecture and SL fibrations}
\label{bp7}

{\it Mirror Symmetry} is a mysterious relationship between pairs 
of Calabi--Yau 3-folds $M,\hat M$, arising from a branch of 
physics known as {\it String Theory}, and leading to some very 
strange and exciting conjectures about Calabi--Yau 3-folds, 
many of which have been proved in special cases.

The {\it SYZ Conjecture} is an attempt to explain Mirror Symmetry
in terms of dual ``fibrations'' $f:M\ra B$ and $\hat f:\hat M\ra B$
of $M,\hat M$ by special Lagrangian 3-folds, including singular
fibers. We give brief introductions to String Theory, Mirror 
Symmetry, and the SYZ Conjecture, and then a short survey of
the state of mathematical research into the SYZ Conjecture, 
biased in favor of the author's own interests.

\subsection{String Theory and Mirror Symmetry}
\label{bp71}

String Theory is a branch of high-energy theoretical physics in 
which particles are modeled not as points but as 1-dimensional 
objects -- ``strings'' --  propagating in some background
space-time $S$. String theorists aim to construct a {\it quantum
theory} of the string's motion. The process of quantization is
extremely complicated, and fraught with mathematical difficulties
that are as yet still poorly understood.

The most popular version of String Theory requires the universe
to be 10-dimensional for this quantization process to work.
Therefore, String Theorists suppose that the space we live in 
looks locally like $S=\R^4\t M$, where $\R^4$ is Minkowski space, 
and $M$ is a compact Riemannian 6-manifold with radius of order 
$10^{-33}$cm, the Planck length. Since the Planck length is so 
small, space then appears to macroscopic observers to be 4-dimensional.

Because of supersymmetry, $M$ has to be a {\it Calabi--Yau $3$-fold}. 
Therefore String Theorists are very interested in Calabi--Yau 3-folds. 
They believe that each Calabi--Yau 3-fold $M$ has a quantization, which 
is a {\it Super Conformal Field Theory} (SCFT), a complicated
mathematical object. Invariants of $M$ such as the Dolbeault 
groups $H^{p,q}(M)$ and the number of holomorphic curves in $M$ 
translate to properties of the SCFT.

However, two entirely different Calabi--Yau 3-folds $M$ and $\hat M$ 
may have the {\it same} SCFT. In this case, there are powerful 
relationships between the invariants of $M$ and of $\hat M$ that
translate to properties of the SCFT. This is the idea behind
{\it Mirror Symmetry} of Calabi--Yau 3-folds.

It turns out that there is a very simple automorphism of the
structure of a SCFT --- changing the sign of a $\U(1)$-action
--- which does {\it not\/} correspond to a classical automorphism 
of Calabi--Yau 3-folds. We say that $M$ and $\hat M$ are {\it mirror}
Calabi--Yau 3-folds if their SCFT's are related by this automorphism.
Then one can argue using String Theory that
\begin{equation*}
H^{1,1}(M)\cong H^{2,1}(\hat M) \quad\text{and}\quad
H^{2,1}(M)\cong H^{1,1}(\hat M).
\end{equation*}
Effectively, the mirror transform exchanges even- and
odd-dimensional cohomology. This is a very surprising result!

More involved String Theory arguments show that, in effect,
the Mirror Transform exchanges things related to the complex
structure of $M$ with things related to the symplectic structure
of $\hat M$, and vice versa. Also, a generating function for the 
number of holomorphic rational curves in $M$ is exchanged
with a simple invariant to do with variation of complex 
structure on $\hat M$, and so on.

Because the quantization process is poorly understood and
not at all rigorous --- it involves non-convergent path-integrals 
over horrible infinite-dimensional spaces --- String Theory
generates only conjectures about Mirror Symmetry, not proofs.
However, many of these conjectures have been verified in
particular cases.

\subsection{Mathematical interpretations of Mirror Symmetry}
\label{bp72}

In the beginning (the 1980's), Mirror Symmetry seemed mathematically 
completely mysterious. But there are now two complementary conjectural
theories, due to Kontsevich and Strominger--Yau--Zaslow, which explain 
Mirror Symmetry in a fairly mathematical way. Probably both are true, 
at some level.

The first proposal was due to Kontsevich \cite{Kont} in 1994. This 
says that for mirror Calabi--Yau 3-folds $M$ and $\hat M$, the derived 
category $D^b(M)$ of coherent sheaves on $M$ is equivalent to the derived 
category $D^b({\rm Fuk}(\hat M))$ of the Fukaya category of $\hat M$, 
and vice versa. Basically, $D^b(M)$ has to do with $M$ as a complex
manifold, and $D^b({\rm Fuk}(\hat M))$ with $\hat M$ as a symplectic 
manifold, and its Lagrangian submanifolds. We shall not discuss this here.

The second proposal, due to Strominger, Yau and Zaslow \cite{SYZ} in 1996, 
is known as the {\it SYZ Conjecture}. Here is an attempt to state it.
\medskip

\noindent{\bf The SYZ Conjecture} {\it Suppose $M$ and\/ $\hat M$ are 
mirror Calabi--Yau $3$-folds. Then (under some additional conditions) 
there should exist a compact topological\/ $3$-manifold\/ $B$ and 
surjective, continuous maps $f:M\ra B$ and\/ $\hat f:\hat M\ra B$, 
such that
\begin{itemize}
\item[{\rm(i)}] There exists a dense open set\/ $B_0\subset B$, such 
that for each\/ $b\in B_0$, the fibers $f^{-1}(b)$ and\/ $\hat f^{-1}(b)$
are nonsingular special Lagrangian $3$-tori $T^3$ in $M$ and\/ $\hat M$.
Furthermore, $f^{-1}(b)$ and\/ $\hat f^{-1}(b)$ are in some sense
dual to one another.
\item[{\rm(ii)}] For each\/ $b\in \De=B\sm B_0$, the fibers $f^{-1}(b)$ 
and\/ $\hat f^{-1}(b)$ are expected to be singular special Lagrangian
$3$-folds in $M$ and\/~$\hat M$.
\end{itemize}}
\medskip

We call $f$ and $\hat f$ {\it special Lagrangian fibrations}, and
the set of singular fibers $\De$ is called the {\it discriminant}.
In part (i), the nonsingular fibers of $f$ and $\hat f$ are supposed
to be {\it dual tori}. What does this mean?

On the topological level, we can define duality between two tori 
$T,\hat T$ to be a choice of isomorphism $H^1(T,\Z)\cong H_1(\hat T,\Z)$. 
We can also define duality between tori equipped with flat Riemannian
metrics. Write $T=V/\La$, where $V$ is a Euclidean vector space
and $\La$ a {\it lattice} in $V\!$. Then the dual torus $\hat T$ is 
defined to be $V^*/\La^*$, where $V^*$ is the dual vector space and 
$\La^*$ the dual lattice. However, there is no notion of duality
between non-flat metrics on dual tori.

Strominger, Yau and Zaslow argue only that their conjecture holds
when $M,\hat M$ are close to the ``large complex structure limit''.
In this case, the diameters of the fibers $f^{-1}(b),\hat f^{-1}(b)$
are expected to be small compared to the diameter of the base space
$B$, and away from singularities of $f,\hat f$, the metrics on the
nonsingular fibers are expected to be approximately flat. 

So, part (i) of the SYZ Conjecture says that for $b\in B\sm B_0$, 
$f^{-1}(b)$ is approximately a flat Riemannian 3-torus, and $\hat f^{-1}(b)$ 
is approximately the dual flat Riemannian torus. Really, the SYZ Conjecture
makes most sense as a statement about the limiting behavior of 
{\it families} of mirror Calabi--Yau 3-folds $M_t$, $\hat M_t$ 
which approach the ``large complex structure limit'' as~$t\ra 0$.

\subsection{The symplectic topological approach to SYZ}
\label{bp73}

The most successful approach to the SYZ Conjecture so far could 
be described as {\it symplectic topological}. In this approach, 
we mostly forget about complex structures, and treat $M,\hat M$ 
just as {\it symplectic manifolds}. We mostly forget about the 
`special' condition, and treat $f,\hat f$ just as {\it Lagrangian 
fibrations}. We also impose the condition that $B$ is a {\it smooth\/} 
3-manifold and $f:M\ra B$ and $\hat f:\hat M\ra B$ are {\it smooth 
maps}. (It is not clear that $f,\hat f$ can in fact be smooth at 
every point, though).

Under these simplifying assumptions, Gross \cite{Gros1,Gros2,Gros3,Gros4}, 
Ruan \cite{Ruan1,Ruan2}, and others have built up a beautiful, detailed 
picture of how dual SYZ fibrations work at the global topological level, 
in particular for examples such as the quintic and its mirror, and for 
Calabi--Yau 3-folds constructed as hypersurfaces in toric 4-folds, using 
combinatorial data.

\subsection{Local geometric approach, and SL singularities}
\label{bp74}

There is also another approach to the SYZ Conjecture, begun by
the author in \cite{Joyc10,Joyc12}, and making use of the ideas
and philosophy set out in Section~\ref{bp6}. We could describe it as
a {\it local geometric} approach.

In it we try to take the special Lagrangian condition seriously 
from the outset, and our focus is on the local behavior of special 
Lagrangian submanifolds, and especially their singularities, rather 
than on global topological questions. Also, we are interested in what 
fibrations of {\it generic} (almost) Calabi--Yau 3-folds might look like.

One of the first-fruits of this approach has been the understanding
that for {\it generic} (almost) Calabi--Yau 3-folds $M$, special
Lagrangian fibrations $f:M\ra B$ will not be smooth maps, but only
piecewise smooth. Furthermore, their behavior at the singular set
is rather different to the smooth Lagrangian fibrations discussed
in Section~\ref{bp73}.

For smooth special Lagrangian fibrations $f:M\ra B$, the discriminant
$\De$ is of codimension 2 in $B$, and the typical singular fiber is 
singular along an ${\cal S}^1$. But in a generic special Lagrangian 
fibration $f:M\ra B$ the discriminant $\De$ is of codimension 1 in $B$, 
and the typical singular fiber is singular at finitely many points.

One can also show that if $M,\hat M$ are a mirror pair of generic
(almost) Calabi--Yau 3-folds and $f:M\ra B$ and $\hat f:\hat M\ra B$ are
dual special Lagrangian fibrations, then in general the discriminants
$\De$ of $f$ and $\hat\De$ of $\hat f$ cannot coincide in $B$, because
they have different topological properties in the neighborhood of a
certain kind of codimension 3 singular fiber. 

This contradicts part (ii) of the SYZ Conjecture, as we have stated it
in Section~\ref{bp72}. In the author's view, these calculations support the 
idea that the SYZ Conjecture in its present form should be viewed 
primarily as a limiting statement, about what happens at the ``large 
complex structure limit'', rather than as simply being about pairs of 
Calabi--Yau 3-folds. A similar conclusion is reached by Mark Gross
in~\cite[\S 5]{Gros4}.

\subsection{$\U(1)$-invariant SL fibrations in $\C^3$}
\label{bp75}

We finish by describing work of the author in \cite[\S 8]{Joyc10}
and \cite{Joyc12}, which aims to describe what the singularities
of SL fibrations of {\it generic} (almost) Calabi--Yau 3-folds
look like, providing they exist.

This proceeds by first studying SL fibrations of subsets of
$\C^3$ invariant under the $\U(1)$-action \eq{bp3eq2}, using
the ideas of Section~\ref{bp35}. For a brief survey of the
main results, see \cite{Joyc11}. Then we argue that the kinds
of singularities we see in codimension 1 and 2 in generic
$\U(1)$-invariant SL fibrations in $\C^3$, also occur in
codimension 1 and 2 in SL fibrations of generic (almost)
Calabi--Yau 3-folds.

Following \cite[Def.~8.1]{Joyc10}, we use the results of
Section~\ref{bp35} to construct a family of SL 3-folds $N_{\bs\al}$ in
$\C^3$, depending on boundary data~$\Phi({\bs\al})$.

\begin{dfn} Let $S$ be a strictly convex domain in $\R^2$ invariant
under $(x,y)\mapsto(x,-y)$, let $U$ be an open set in $\R^3$, and
$\al\in(0,1)$. Suppose $\Phi:U\ra C^{3,\al}(\pd S)$ is a continuous
map such that if $(a,b,c)\ne(a,b',c')$ in $U$ then $\Phi(a,b,c)-
\Phi(a,b',c')$ has exactly one local maximum and one local minimum
in~$\pd S$.

For ${\bs\al}=(a,b,c)\in U$, let $f_{\bs\al}\in C^{3,\al}(S)$
or $C^1(S)$ be the unique (weak) solution of \eq{bp3eq6} with
$f_{\bs\al}\vert_{\pd S}=\Phi({\bs\al})$, which exists by
Theorem \ref{bp3thm1}. Define
\begin{equation*}
u_{\bs\al}=\frac{\pd f_{\bs\al}}{\pd y}
\quad\hbox{and}\quad
v_{\bs\al}=\frac{\pd f_{\bs\al}}{\pd x}.
\end{equation*}
Then $(u_{\bs\al},v_{\bs\al})$ is a solution of \eq{bp3eq5} in
$C^{2,\al}(S)$ if $a\ne 0$, and a weak solution of \eq{bp3eq4}
in $C^0(S)$ if $a=0$. Also $u_{\bs\al},v_{\bs\al}$ depend
continuously on ${\bs\al}\in U$ in $C^0(S)$, by Theorem~\ref{bp3thm1}.

For each ${\bs\al}=(a,b,c)$ in $U$, define $N_{\bs\al}$ in $\C^3$ by
\e
\begin{split}
N_{\bs\al}=\bigl\{(z_1,z_2,z_3)\in\C^3:\,&
z_1z_2=v_{\bs\al}(x,y)+iy,\quad z_3=x+iu_{\bs\al}(x,y),\\
&\ms{z_1}-\ms{z_2}=2a,\quad (x,y)\in S^\circ\bigr\}.
\end{split}
\label{bp7eq1}
\e
Then $N_{\bs\al}$ is a noncompact SL 3-fold without boundary in $\C^3$,
\label{bp7def}
which is nonsingular if $a\ne 0$, by Proposition~\ref{bp3prop1}.
\end{dfn}

In \cite[Th.~8.2]{Joyc10} we show that the $N_{\bs\al}$ are the
fibers of an {\it SL fibration}.

\begin{thm} In the situation of Definition \ref{bp7def}, if\/
${\bs\al}\ne{\bs\al}'$ in $U$ then $N_{\bs\al}\cap N_{{\bs\al}'}
=\emptyset$. There exists an open set\/ $V\subset\C^3$ and a
continuous, surjective map $F:V\ra U$ such that\/ $F^{-1}({\bs\al})
=N_{\bs\al}$ for all\/ ${\bs\al}\in U$. Thus, $F$ is a special
Lagrangian fibration of\/ $V\subset\C^3$, which may include
\label{bp7thm}
singular fibers.
\end{thm}

It is easy to produce families $\Phi$ satisfying Definition
\ref{bp7def}. For example \cite[Ex.~8.3]{Joyc10}, given any
$\phi\in C^{3,\al}(\pd S)$ we may define $U=\R^3$ and $\Phi:
\R^3\ra C^{3,\al}(\pd S)$ by $\Phi(a,b,c)=\phi+bx+cy$. So this
construction produces very large families of $\U(1)$-invariant
SL fibrations, including singular fibers, which can have any
multiplicity and type.

Here is a simple, explicit example. Define $F:\C^3\ra\R\t\C$ by
\e
\begin{gathered}
F(z_1,z_2,z_3)=(a,b),\quad\text{where}\quad 2a=\ms{z_1}-\ms{z_2} \\
\text{and}\quad
b=\begin{cases}
z_3, & a=z_1=z_2=0, \\
z_3+\bar z_1\bar z_2/\md{z_1}, & a\ge 0,\; z_1\ne 0,\\
z_3+\bar z_1\bar z_2/\md{z_2}, & a<0.
\end{cases}
\end{gathered}
\label{bp7eq2}
\e
This is a piecewise-smooth SL fibration of $\C^3$. It is not
smooth on~$\md{z_1}=\md{z_2}$.

The fibers $F^{-1}(a,b)$ are $T^2$-cones singular at $(0,0,b)$ 
when $a=0$, and nonsingular ${\cal S}^1\t\R^2$ when $a\ne 0$. 
They are isomorphic to the SL 3-folds of Example \ref{bp3ex1}
under transformations of $\C^3$, but they are assembled to
make a fibration in a novel way.

As $a$ goes from positive to negative the fibers undergo a surgery, 
a Dehn twist on ${\cal S}^1$. The reason why the fibration is only 
piecewise-smooth, rather than smooth, is really this topological 
transition, rather than the singularities themselves. The fibration 
is not differentiable at every point of a singular fiber, rather
than just at singular points, and this is because we are jumping
from one moduli space of SL 3-folds to another at the singular fibers.

I claim that $F$ is a local model for codimension one 
singularities of SL fibrations of generic almost Calabi--Yau
3-folds. The reason for this is that these $T^2$-cone
singularities are {\it stable}, as in Definition \ref{bp6def3},
so SL 3-folds $X$ with these singularities form {\it smooth
moduli spaces} $\M_\sX$ by Corollary~\ref{bp6cor1}.

The singularities are automatically {\it transverse},
as in Definition \ref{bp6def6}, so we can apply
\cite[Th.~8.10]{Joyc17} to compute the {\it index}
$\ind(X)$ of the singularities, as in Section \ref{bp65}.
This is done in detail in \cite[\S 10]{Joyc17}.
If the topology of $X$ is suitably chosen then
$\ind(X)=1$, so $\M_\sX$ has codimension one in $\M_\sN$.
The singular behavior is stable under small exact
perturbations of the underlying almost Calabi--Yau structure.

I also have a $\U(1)$-invariant model for codimension two 
singularities, described in \cite{Joyc12}, in which two of the
codimension one $T^2$-cones come together and cancel out. I
conjecture that it too is a typical codimension two singular
behavior in SL fibrations of generic almost Calabi--Yau 3-folds.
I do not expect codimension three singularities in generic SL
fibrations to be locally $\U(1)$-invariant, and so this
approach will not help.

\end{document}